\documentclass[a4paper,10pt]{amsart}
\def\C{\mathbb{C}}

\def\R{\mathbb{R}}

\def\B{\mathcal B}

\def\mb{\mathbb}

\def\beqa{\begin{eqnarray*}}
\def\eeqa{\end{eqnarray*}}

 \newtheorem{thm}{Theorem}
 \newtheorem{cor}[thm]{Corollary}
 \newtheorem{lem}[thm]{Lemma}
 \newtheorem{prop}[thm]{Proposition}
 \newtheorem{defn}[thm]{Definition}
 
 \newtheorem{ex}[thm]{Example}
\newcommand{\be}{\begin{equation}}
\newcommand{\ee}{\end{equation}}
\newcommand{\bea}{\begin{eqnarray}}

\newcommand{\eea}{\end{eqnarray}}
\newcommand{\Bea}{\begin{eqnarray*}}
\newcommand{\Eea}{\end{eqnarray*}}

\newcounter{cnt1}
\newcounter{cnt2}
\newcounter{cnt3}
\newcommand{\blr}{\begin{list}{$($\roman{cnt1}$)$}
 {\usecounter{cnt1} \setlength{\topsep}{0pt}
 \setlength{\itemsep}{0pt}}}
\newcommand{\bla}{\begin{list}{$($\alph{cnt2}$)$}
 {\usecounter{cnt2} \setlength{\topsep}{0pt}
 \setlength{\itemsep}{0pt}}}
\newcommand{\bln}{\begin{list}{$($\arabic{cnt3}$)$}
 {\usecounter{cnt3} \setlength{\topsep}{0pt}
 \setlength{\itemsep}{0pt}}}
\newcommand{\el}{\end{list}}

\date{}
\begin{document}

\author{B. V. Rajarama Bhat and  Mithun  Mukherjee}

\title{ Inclusion systems and amalgamated products of product
systems }
\maketitle
\begin{center}
Indian Statistical Institute,\\
R. V. College Post,\\
Bangalore-560059, India.\\
bhat@isibang.ac.in, mithun@isibang.ac.in\\

\vskip 10pt

{\bf Abstract\footnote {To appear in the Journal  `Infinite Dimensional Analysis, Quantum Probability and 
Related Topics' (2010).}}

\end{center}

\vskip 4pt


Here we generalize the concept of Skeide product, introduced by
Skeide, of two product systems via a pair of normalized units. This
new notion is called amalgamated  product of product systems, and
now the amalgamation can be done using  contractive morphisms. Index
of amalgamation product (when done through units) adds up for
normalized units but for non-normalized units, the index is one more
than the sum. We define inclusion systems and use it as a tool for
index computations. It is expected that this notion will have other
uses.

\section{Introduction}

Studying quantum dynamical semigroups (or completely positive
semigroups) and their dilations is important in understanding
irreversible quantum dynamics. In this context, R.T. Powers posed
the following problem at the 2002 AMS summer conference on `Advances
in Quantum Dynamics' held at Mount Holyoke: Let $\mathcal {B}
(\mathcal {H})$ and $\mathcal {B}(\mathcal {K})$ be algebras of all
bounded operators on two Hilbert spaces $\mathcal H$ and $\mathcal
K.$ Suppose $\phi =\{ \phi _t: t\geq 0\} $ and $\psi = \{ \psi _t:
t\geq 0\} $ are two contractive completely positive (CP) semigroups
on $\mathcal {B} (\mathcal {H})$ and $\mathcal {B}(\mathcal {K})$
respectively and $U=\{ U_t: t\geq 0\}$ and $V=\{ V_t: t\geq 0\}$ are
two strongly continuous semigroups of isometries which intertwine
$\phi _t$ and $\psi _t$ respectively. Consider the CP semigroup
$\tau _t$ on $\mathcal {B}( \mathcal {H}\oplus \mathcal {K})$
defined by
 $\tau _t\left(
 \begin{array}{cc}
 X & Y \\
 Z & W \\
 \end{array}
 \right)=\left(
 \begin{array}{cc}
 \phi _t(X) & U_tYV^*_t \\
 V_tZU^*_t & \psi _t(W) \\
 \end{array}
 \right).$
 How is the minimal dilation (in the sense of  \cite{Bh01},
\cite{Bh02} ) of $\tau  $ related to minimal dilations  of $\phi $
and $\psi $? In fact, Powers was interested in a more specific
question. It is the following.  Recall that by W. Arveson
\cite{Ar01} we can associate a tensor product system of Hilbert
spaces with every $E$-semigroup of ${\mathcal B} ({\mathcal K}).$
Since the minimal dilation is unique we might say that we are
associating a product system  to a given contractive CP semigroup
\cite{Bh01}. (This can also be done more directly as in \cite{BS}).
The question was to `What is the product system of $\tau $ in terms
of product systems of $\phi $ and $\psi .$ Is it the tensor
product?' This was answered by Skeide (See \cite{BLS}, \cite{Po},
\cite{Sk02}, \cite{Sk03}). It turns out that though the index of product system
of $\tau  $ is sum of indices of $\phi $ and $\psi $,  it is not
tensor product but a product introduced by Skeide in \cite{Sk01}. We
will be calling it as Skeide product.   It is a kind of amalgamated
tensor product where we identify two units. Strangely proof of this
fact depended upon the fact that units (intertwining semigroups)
$\{U_t \}, \{V_t\}$ are normalized, though $\tau  $ can be
constructed even when they are just contractive. This raised the
question as to what is the product system when the units are not
normalized.  We answer this question and we see that surprisingly
the index increases by 1 if the units are not normalized. The
original motivation in \cite{Sk01} was to get an appropriate product
for product systems Hilbert modules where the index is additive,
since there is no obvious notion of tensor products for product
systems of Hilbert modules. Since this article is only about Hilbert
spaces, we do not dwell more on this point.

To begin with we introduce the notion of inclusion systems. These
are parametrized families of Hilbert spaces exactly like product
systems except that now unitaries are replaced by isometries.
Actually these objects seem to be ubiquitous in the field of product
systems. Even while associating product systems to CP semigroups
what one gets first are inclusion systems, and then an inductive
limit procedure gives product system \cite{BS}. In \cite{BS},
\cite{BBLS} this procedure has been elaborated and has been
exploited, in the more general context of product systems of Hilbert
modules. Since the same construction is being repeatedly used, it is
good to extract the essence of the method and put it in a general
framework. This is what is being done here. In other words we define
inclusion systems and show that every inclusion system gives rise to
a product system in a natural way (by taking inductive limits). It
is remarkable here that basic properties of the product systems such
as existence/non-existence of units, structure of morphisms etc. can
be read of at the level of inclusion systems. This is the
observation which we wish to stress.  We believe that this technique
will be very useful in many other contexts. While writing this
article, we have come to know that Shalit and Solel \cite{SS}, have
called inclusion systems as `subproduct systems' and have looked at
their general theory, connections with CP semigroups, even
subproduct systems of correspondences.  As per authors of \cite{SS},
their goal is to extend the dilation theory to more general
semigroups, that is, instead of $\mathbb Z_+, \mathbb R_+$ they wish
to consider $\mathbb Z_+^d, \mathbb R_+^d$ etc. for the time
parameter. They restrict mostly to discrete semigroups as  the
multivariable theory is quite involved already for discrete
semigroups. Consequently inductive limits do not appear in their
work. So the results in \cite{SS} are quite different from ours.
Since there is the possibility of confusion between `subproduct
systems' and `product subsystems', we avoid this terminology.
Occasionally we need to talk about one inclusion system contained in
another one, and here too the usage `subproduct system' becomes a
bit awkward. Our work and some questions arose from it have
motivated Tsirelson to look deeper into subproduct systems with
finite dimensional spaces (\cite{Tsi01}, \cite{Tsi02}).

In Section 3 we describe 'amalgamation' of two inclusion systems via
a contractive morphism. This is more general than Skeide product and
is useful for treating the case of more general `corners' in Powers'
problem. There is an interesting relationship between units of the
amalgamated product and units of the individual systems (Lemma 19).
Lastly we treat the special case, where the amalgamation is done
through a pair of units, using the contraction morphism
$D_t=|u^0_t\rangle \langle v^0_t|$ where $u^0$ and $v^0$ are two
fixed units of product systems  $\mathcal E$ and $\mathcal F$ such
that $\parallel u^0_t
\parallel
\parallel v^0_t \parallel \leq 1$ for all $t>0$. Then it turns out that $Ind(\mathcal E \otimes \mathcal F)=Ind(\mathcal E)+ Ind(\mathcal
F)$ if both the units are normalized, $Ind(\mathcal E \otimes _D
\mathcal F)=Ind(\mathcal E) + Ind(\mathcal F) + 1$ otherwise.

\section{Inclusion Systems}\label{inclusion}

\begin{defn}
An Inclusion System $(E, \beta )$ is  a family of Hilbert spaces $E=
\{ E_t,t\in (0, \infty )\}$ together with isometries
$\beta_{s,t}$:$E_{s+t}\rightarrow E_s\otimes E_t$, for   $s,t \in
(0, \infty )$, such that $\forall$ $r,s,t \in (0, \infty),$ ~
$(\beta_{r,s}\otimes 1_{E_t})\beta_{r+s,t}=(1_{E_r}\otimes
\beta_{s,t})\beta_{r,s+t}.$ It is said to be a product system if
further every $\beta _{s, t}$ is a unitary.
\end{defn}
At the moment we are not putting any measurability conditions on
inclusion systems. Such technical conditions can be put when they
become necessary. Of course, every product system is an inclusion
system. Now here, there is a subtle point to note. Defining
unitaries for product systems usually go from $E_s\otimes E_t$ to
$E_{s+t}$ and are associative. We have taken their adjoint maps
which are `co-associative'. So one might say that we are actually
looking at `co-product systems' and abusing the terminology by
calling them `product systems'. We thank the referee for pointing
this out to us.

Here are some genuine inclusion systems.

\begin{ex}
Take $E_t\equiv \mathbb C^2$ with ortho-normal basis $\{ e_0,
e_1\}$. Define $\beta _{s,t}:E_{s+t}\to E_s\otimes E_t$ by
$$\beta _{s,t}e_0=e_0\otimes e_0,~~ \beta _{s,t}(e_1)=
\frac{1}{\sqrt{s+t}}(\sqrt se_1\otimes e_0+\sqrt{t} e_0\otimes
e_1).$$ Then $(E, \beta )$ is an inclusion system.
\end{ex}

In the following example we are making use of concepts such as `units' and spatial product systems (those which have units) as in \cite{Ar03}. The reader may also refer to Definition 8, below.

\begin{ex}
Let $({\mathcal  F}, C)$ be a spatial product system.  Let
${\mathcal U }^{\mathcal F}$ be the set of units of this product
system. Now $(E, \beta )$ with $E_t= \overline {\mbox {span}} \{u_t:
u\in {\mathcal U }^{\mathcal F}\}$, and $\beta _{s, t}= C_{s,
t}|_{E_{s+t}}$ is an inclusion system.
\end{ex}

Stinespring dilations of semigroups of completely positive maps is another source of inclusion systems. This
will be explained towards the end of this Section.

Our first job is to show that every inclusion system leads to a product system in a natural way. Here we explain
this procedure. So consider an inclusion system $(E, \beta )$. Let for $t \in \R_+, $ $J_t=\{(t_n,t_{n-1},\ldots
, t_1): t_i>0, \sum _{i=1}^n t_i=t, n\geq 1\}$. \ For $\textbf{s}=(s_m,s_{m-1},\ldots , s_1) \in J_s$, and
$\textbf{t}=(t_n,t_{n-1}, \ldots , t_1) \in J_t$ we define $\textbf{s} \smile \textbf{t}:=(s_m,s_{m-1},\ldots
,s_1,t_n,t_{n-1}, \ldots , t_1) \in J_{s+t}$. Now fix $t \in \R_+$. On $J_t$ define a partial order
$\textbf{t}\geq\textbf{s}=(s_m,s_{m-1},\ldots, s_1)$ if for each $i,$ $(1\leq i\leq m)$ there exists (unique)
$\textbf{s}_i \in J_{s_i}$ such that $\textbf{t}=\textbf{s}_m\smile\textbf{s}_{m-1}\smile \cdots
\smile\textbf{s}_1$. For $\textbf{t}= (t_n, t_{n-1}, \ldots t_1) $ in $J_t$ define $E_\textbf{t}=E_{t_n}\otimes
E_{t_{n-1}}\otimes \cdots \otimes E_{t_1}$. For $\textbf {s}=(s_m,\ldots
,s_1)\leq\textbf{t}=(\textbf{s}_m\smile\cdots \smile\textbf{s}_1)$ in $J_t,$ define
$\beta_{\textbf{t},\textbf{s}}:E_\textbf{s}\rightarrow E_\textbf{t} $ by
$\beta_{\textbf{t},\textbf{s}}=\beta_{\textbf{s}_m,s_m}\otimes\beta_{\textbf{s}_{m-1},s_{m-1}}\otimes \cdots
\otimes\beta_{\textbf{s}_1,s_1}$ where we define $\beta_{\textbf{s},s}:E_s\rightarrow E_\textbf{s}$ inductively
as follows: Set $\beta_{s,s}=id_{E_s}$. For $\textbf{s}=(s_m,s_{m-1},\ldots ,s_1)$, $\beta_{\textbf{s},s}$ is
the composition of maps:
$$
 (\beta _{s_m, s_{m-1}}\otimes I)(\beta _{s_m+s_{m-1}, s_{m-2}}\otimes I)\cdots (\beta _{s_{m}+\cdots +s_3, s_2}\otimes I)\beta _{s_m+\cdots +s_2, s_1}.
$$

\begin{lem} Let $t\in \R_+$ be
fixed and consider the partially ordered set $J_t$ defined above. Then
$\{E_\textbf{t},\beta_{\textbf{s},\textbf{r}}: \textbf{r},\textbf{s},\textbf{t} \in J_t $\} forms an Inductive
System of Hilbert spaces in the  sense that: (i)  $\beta_{\textbf{s},\textbf{s}}=id_{E_\textbf{s}}$ for
$\textbf{s}\in J_t$; (ii)
$\beta_{\textbf{t},\textbf{s}}\beta_{\textbf{s},\textbf{r}}=\beta_{\textbf{t},\textbf{r}}$ for
$\textbf{r}\leq\textbf{s}\leq\textbf{t}\in J_t.$
\end{lem}
Proof: Only (ii) needs to be proved.  Let $\textbf{r}=(r_n, \ldots , r_1)$,
 $\textbf{s}=\textbf{r}_n\smile \cdots \smile\textbf{r}_1$, where $\textbf{r}_i=(r_{ik_i}, \ldots ,r_{i1}),$ $1\leq i\leq n.$
 So
$$\textbf{t}=(\textbf{r}_{nk_n}\smile \ldots \smile \textbf{r}_{n1})\smile(\textbf{r}_{(n-1)k_{n-1}}\smile \ldots
\smile\textbf{r}_{(n-1)1})\smile \ldots \smile(\textbf{r}_{1k_1}\smile \ldots \smile  \textbf{r}_{11}).$$ Now
\Bea \beta_{\textbf{t},\textbf{s}}\beta_{\textbf{s},\textbf{r}}&=&
\beta_{\textbf{t},\textbf{s}}(\beta_{\textbf{r}_n,r_n}\otimes \cdots \otimes\beta_{\textbf{r}_1,r_1})\\ &=&
(\beta_{\textbf{r}_{nk_n},r_{nk_n}}\otimes\beta_{\textbf{r}_{nk_{n-1}},r_{nk_{n-1}}}\otimes \cdots
\otimes\beta_{\textbf{r}_{n1},r_{n1}}\otimes \cdots \otimes\beta_{\textbf{r}_{1k_1},r_{1k_1}}\otimes \cdots
\otimes\beta_{\textbf{r}_{11},r_{11}})\\ &&(\beta_{\textbf{r}_n,r_n}\otimes \cdots \otimes\beta_{\textbf{r}_1,r_1})\\
&=& [\beta_{(\textbf{r}_{nk_n}\smile \cdots
\smile\textbf{r}_{n1}),\textbf{r}_n}\otimes\beta_{(\textbf{r}_{(n-1)k_{n-1}}\smile \cdots
\smile\textbf{r}_{(n-1)1}),\textbf{r}_{n-1}}\otimes \cdots \otimes\beta_{(\textbf{r}_{1k_1}\smile \cdots
\smile\textbf{r}_{11}),\textbf{r}_1}]\\ &&(\beta_{\textbf{r}_n,r_n}\otimes \cdots \otimes\beta_{\textbf{r}_1,r_1})\\
&=& \beta_{(\textbf{r}_{nk_n}\smile \cdots
\smile\textbf{r}_{n1}),r_n}\otimes\beta_{(\textbf{r}_{(n-1)k_{n-1}}\smile \cdots
\smile\textbf{r}_{(n-1)1}),r_{n-1}}\otimes\cdots \beta_{(\textbf{r}_{1k_1}\smile \cdots
\smile\textbf{r}_{11}),r_1}\\ &=& \beta_{\textbf{t},\textbf{r}} \Eea \qed

\begin{thm}\label{inductive }
Suppose $(E, \beta )$ is an inclusion system. Let $\mathcal E_t=indlim_{J_t} E_\textbf{s}$ be the inductive
limit of $E_\textbf{s}$ over $J_t$ for $t>0.$ Then $E= \{\mathcal E_t : t>0\}$ has the structure of a product
system of Hilbert spaces.
\end{thm}

Proof: We recall four basic properties of the inductive limit construction.  (i) There exist canonical
isometries $i_\textbf{s}:E_\textbf{s}\rightarrow \mathcal E_t$  such that given $\textbf{r}$ , $\textbf{s}$ in
$J_t$ with $\textbf{r}\leq \textbf{s}$ , $i_\textbf{s}\beta_{\textbf{s},\textbf{r}}=i_\textbf{r}$. ~(ii)
$\overline {\mbox{span}} \{i_\textbf{s}(a):a \in E_\textbf{s}, \textbf{s}\in J_t\}=\mathcal E_t.$ (iii) The
following universal property holds : Given a Hilbert space $\mathcal G$ and isometries $g_\textbf{s}:
E_\textbf{s}\rightarrow \mathcal G$  satisfying consistency condition
$g_\textbf{s}\beta_{\textbf{s},\textbf{r}}=g_\textbf{r}$ for all $\textbf{r}\leq \textbf{s} $ there exists a
unique isometry $g:\mathcal E_t \rightarrow \mathcal G$ such that $g_\textbf{s}=gi_\textbf{s}$ $\forall
\textbf{s} \in J_t$. (iv) Suppose $K \subseteq J_t$ has the following property: Given $\textbf{s} \in J_t$ there
exists $\textbf{t} \in K$ such that $\textbf{s} \leq \textbf{t}$. Then $K$ is indeed a directed set with the
order inherited from $J_t$ and that the injection $K\rightarrow J_t$ is a cofinal function. In other words,
$(x_{\textbf{s}})_{\textbf{s}\in K}$ is a subnet of $(x_{\textbf{t}})_{\textbf{t}\in J_t}.$
$\mbox{indlim}_{J_t}E_\textbf{s}=\mbox{indlim}_KE_\textbf{s}$.

Define $J_s\smile J_t=\{\textbf{s}\smile\textbf{t}:\textbf{s}\in J_s ,\textbf{t} \in J_t \}$. Given any element
$\textbf{r} \in J_{s+t},$ there exist  $\textbf{s} \in J_s$ and $\textbf{t} \in J_t$ such that $\textbf{s}\smile
\textbf{t} \geq \textbf{r}.$ So by the second property quoted above $\mathcal E_{s+t}= \mbox{indlim}_{J_s\smile
J_t}E_{\textbf{s}\smile \textbf{t}}=\mbox{indlim}_{\textbf{s}\smile \textbf{t}\in J_s\smile
J_t}E_\textbf{s}\otimes E_\textbf{t}.$ Let $i_\textbf{s}: E_\textbf{s}\rightarrow \mathcal E_s$ , $i_\textbf{t}:
E_\textbf{t}\rightarrow \mathcal E_t$ be the canonical isometries. Consider the map $i_\textbf{s}\otimes
i_\textbf{t}:E_{\textbf{s}\smile \textbf{t}}\rightarrow \mathcal E_s\otimes \mathcal E_t$ , for
$\textbf{s}\smile \textbf{t}\in J_s\smile J_t$. Note that $\textbf{s}^\prime \smile \textbf{t}^\prime \leq
\textbf{s} \smile \textbf{t}$ in $J_s \smile J_t$
 implies $\textbf{s}^\prime \leq \textbf{s},$
$\textbf{t}^\prime \leq \textbf{t}.$ Now as $\beta_{\textbf{s}\smile \textbf{t},\textbf{s}^\prime\smile
\textbf{t}^\prime}=\beta_{\textbf{s},\textbf{s}^\prime}\otimes \beta_{\textbf{t},\textbf{t}^\prime},$ we get
$(i_\textbf{s}\otimes i_\textbf{t})\beta_{\textbf{s}\smile\textbf{t},\textbf{s}^\prime\smile\textbf{t}^\prime} =
i_\textbf{s}\beta_{\textbf{s},\textbf{s}^\prime}\otimes i_{\textbf{t}}\beta_{\textbf{t},\textbf{t}^\prime}
=i_{\textbf{s}^\prime}\otimes i_{\textbf{t}^\prime}$ So by the universal property, there is a unique isometry
$B_{s,t}:\mathcal E_{s+t}\rightarrow \mathcal E_s\otimes \mathcal E_t$ such that
$B_{s,t}i_{\textbf{s}\smile\textbf{t}}=i_\textbf{s}\otimes i_\textbf{t}.$ From (ii) it is clear that $B_{s,t}$
is a unitary map from $\mathcal E_{s+t}$ to $\mathcal E_s\otimes \mathcal E_t.$

Now to check $(B_{r,s}\otimes 1_{\mathcal E_t})B_{r+s,t}=(1_{\mathcal E_r}\otimes B_{s,t})B_{r,s+t},$ enough to
check it on the vectors of the form $i_{\textbf{r}\smile\textbf{s}\smile\textbf{t}}(a\otimes b\otimes c)$ , $a
\in E_\textbf{r}$, $b \in E_\textbf{s}$ , $c \in E_\textbf{t}$. We have \Bea (B_{r,s}\otimes 1_{\mathcal
E_t})B_{r+s,t}(i_{\textbf{r}\smile\textbf{s}\smile\textbf{t}}(a\otimes b\otimes c)&=& (B_{r,s}\otimes
1_{\mathcal E_t})(i_{\textbf{r}\smile \textbf{s}}(a \otimes b)\otimes
i_\textbf{t}(c))\\
&=& B_{r,s}i_{\textbf{r}\smile \textbf{s}}(a \otimes b)\otimes
i_\textbf{t}(c)\\
&=& i_\textbf{r}(a)\otimes i_\textbf{s}(b)\otimes
i_\textbf{t}(c) \Eea
 And also
\Bea (1_{\mathcal E_r}\otimes B_{s,t})B_{r,s+t}i_{\textbf{r}\smile\textbf{s}\smile\textbf{t}}(a\otimes b\otimes
c)& = &(1_{\mathcal E_r}\otimes B_{s,t})(i_\textbf{r}(a)\otimes i_{\textbf{s}\smile \textbf{t}}(b
\otimes c))\\
& =& i_\textbf{r}(a)\otimes i_\textbf{s}(b)\otimes
i_\textbf{t}(c).\Eea
 This proves the Theorem. \qed

\begin{defn}
Given  an inclusion system $(E, \beta )$, the product system $({\mathcal E}, B)$ constructed as in the previous theorem is called the product system {\em generated by} the inclusion system
$(E, \beta ).$
\end{defn}

It is to be noted that if $(E, \beta )$ is already a product system then its generated system is itself.

\begin{defn}
Let $(E, \beta )$ and $(F, \gamma)$ be two inclusion systems.  Let
$A=\{ A_t: t>0\}$ be a family of linear maps $A_t:E_t\to F_t$,
satisfying $\|A_t\| \leq e^{tk}$ for some $k\in \R.$ Then $A$ is
said to be a {\em morphism} or a {\em weak morphism} from $(E, \beta
)$ to $(F, \gamma )$ if
$$A_{s+t}= \gamma _{s, t}^*(A_s\otimes A_t)\beta _{s, t} ~~\forall s, t>0.$$
It is said to be a {\em strong morphism} if
$$\gamma _{s,t} A_{s, t}= (A_s\otimes A_t)\beta _{s, t} ~~\forall s, t>0.$$
\end{defn}

It is clear that every strong morphism is a weak morphism but the
converse is not true. However the two notions coincide for product
systems, as the linking maps are all unitaries. The exponential
boundedness condition becomes important when we take inductive
limits.  We also note that adjoint of a weak morphism is a weak
morphism. The adjoint of a strong morphism need not be a strong
morphism, but it is at least a weak morphism. Compositions of strong
morphisms is a strong morphism, but this need not be true for weak
morphisms.

\begin{defn}
Let $(E, \beta )$ be an inclusion system. Let $u=\{ u_t: t>0\}$ be a family of vectors such that $(1)$ for all
$t>0,$ $u_t\in E_t$ $(2)$ there is a $k\in \R,$ such that $\|u_t\|\leq e^{tk},$ for all $t>0.$ and  $(3)$
$u_t\neq 0$ for some $t>0.$ Then $u$ is said to be a {\em unit\/} or a {\em weak unit\/} if
$$u_{s+t}= \beta _{s, t}^*(u_s\otimes u_t) ~~\forall s, t>0.$$
It is said to be a {\em strong unit\/} if
$$\beta _{s, t}u_{s+t}=u_s\otimes u_t
 ~~\forall s, t>0.$$
\end{defn}

 A weak (resp. strong) unit of $(E, \beta )$ can be thought of
as a non-zero weak (resp. strong) morphism from the trivial product system $(F, \gamma)$, where $F_t\equiv \C$
and $\gamma _{s, t}(a)= a\otimes 1.$ As any morphism $A:(F,\gamma)\rightarrow (E,\beta)$ is completely
determined by the values $A_t(1),$ $t>0.$ It is easy to see that $(A_t(1))_{t>0}$ is a weak or strong unit
if $A$ is weak or strong morphism respectively.

\begin{lem}
Let $A_t:(E_t , \beta_{r,s})\rightarrow (F_t ,
\beta^\prime_{r,s})$ be a morphism and let $v=(v_t)_{t>0}$ be a unit
of $F$. Then $(A^*_tv_t)_{t>0}$ is a unit, provided $A^*_tv_t\neq 0$ for some $t>0$.
\end{lem}
Proof:  Suppose $\|A_t\|\leq e^{kt}$ and $\|v_t\|\leq e^{lt},$ for
some $l,k\in \R$ Now $$\|A^*_tv_t\| \leq \| A^*_t\| \|v_t\|\leq
e^{(k+l)t}.$$ So $(A^*v)_{t>0} $ is exponentially bounded. Also \Bea
\beta^*_{s,t}(A^*_sv_s\otimes A^*_tv_t)&=&
\beta^*_{s,t}(A^*_s\otimes A^*_t)(v_s\otimes v_t)\\&=&[(A_s\otimes
A_t)\beta_{s,t}]^*(v_s\otimes
v_t)\\&=&[\beta^\prime_{s,t}A_{s+t}]^*(v_s\otimes v_t)\\&=&
A^*_{s+t}{\beta^{\prime*}}_{s,t}(v_s\otimes v_t)\\ &=&
A^*_{s+t}v_{s+t}. \Eea \qed

\begin{thm}\label{unit}
Let $(E , \beta)$ be an inclusion system and let $(\mathcal E , B)$ be the product system generated by it. Then
the canonical map $(i_t)_{t>0}:E_t \rightarrow \mathcal E_t$ is an isometric strong morphism of inclusion
systems. Further $i^*$ is an isomorphism between units of $(\mathcal E, B)$ and units of $(E, \beta ).$

\end{thm}
Proof: The first statement is clear from the construction. For $\textbf{s}=\{s_n,s_{n-1},\cdots,s_1\}\in J_t,$
denote $\mathcal E_\textbf{s}=\mathcal E_{s_n}\otimes\cdots\otimes \mathcal E_{s_1}.$  Let
$i_\textbf{s}:E_\textbf{s}\rightarrow \mathcal E_t$ be the canonical map and let $B_{\textbf{s},t}:\mathcal
E_t\rightarrow\mathcal E_\textbf{s}$ be the unitary map as defined earlier. From the proof of Theorem \ref{inductive },
$B_{\textbf{s},t}i_\textbf{s}=i_{s_n}\otimes i_{s_{n-1}}\cdots\otimes i_{s_1}.$ For any unit $v$ in $(\mathcal
E,B),$ denote $v_\textbf{s}=v_{s_n}\otimes\cdots\otimes v_{s_1}\in \mathcal E_\textbf{s}.$ Then from definition
$B^*_{\textbf{s},t}v_\textbf{s}=v_t.$

Now we prove the injectivity of $i^*.$ Consider two units $v, w$ of $(\mathcal E, B)$ such
 that $i^*_tv_t=i^*_tw_t$  for all $t>0$. Now for any $\textbf{s}\in J_t,$
\Bea i^*_\textbf{s}v_t &=& i^*_\textbf{s}B^*_{\textbf{s},t}v_\textbf{s}\\
&=& (B_{\textbf{s},t}i_\textbf{s})^*v_\textbf{s}\\ &=& (i^*_{s_n}\otimes \cdots \otimes
i^*_{s_1})(v_{s_n}\otimes \cdots \otimes v_{s_1})\\ &=& i^*_{s_n}v_{s_n}\otimes\cdots\otimes i^*_{s_1}v_{s_1}\\
&=& i^*_{s_n}w_{s_n}\otimes\cdots\otimes i^*_{s_1}w_{s_1}\\ &=&
i^*_\textbf{s}w_t.  \Eea This implies
$i_\textbf{s}i^*_\textbf{s}v_t=i_\textbf{s}i^*_\textbf{s}w_t$ for
all $\textbf{s}\in J_t$. Note $\mathcal E_t=\overline {\mbox
{span}}\{i_\textbf{s}(a): a \in E_\textbf{s}, \textbf{s} \in J_t\}.$
The identity
$i_\textbf{t}\beta_{\textbf{t},\textbf{s}}=i_\textbf{s}$ implies
that for $\textbf{s}\leq\textbf{t}\in J_t,$
$$i_\textbf{s}i^*_\textbf{s}i_\textbf{t}i^*_\textbf{t}=i_\textbf{s}i^*_\textbf{s}.$$
So the net of projections $\{i_\textbf{s}i^*_\textbf{s}: \textbf{s}
\in J_t \}$ converges strongly to identity. So we get $v_t=w_t.$ From previous lemma $i^*$ sends  unit to unit, provided the image is non-trivial, which is guaranteed by injectivity of $i^*.$

To see surjectivity, consider a unit  $u$ of the inclusion system
$(E , \beta)$ with $\|u_t\|\leq e^{kt}$ for some $k\in \R_+.$ Fix
$t>0.$ Define $u_\textbf{s}:=u_{s_n}\otimes \cdots \otimes u_{s_1}$
for $(s_n,s_{n-1}, \ldots , s_1) \in J_t.$ Now it follows easily
that for $\textbf{s}\leq \textbf{t},$
$u_\textbf{s}=\beta^*_{\textbf{t},\textbf{s}}u_\textbf{t}$. Now
Consider the  bounded net $\{i_\textbf{s}u_\textbf{s}, \textbf{s}
\in J_t\}.$ For $\textbf{s}\leq\textbf{t}\in J_t,$ we have \Bea
i_\textbf{s}i^*_\textbf{s}i_\textbf{t}u_\textbf{t}&=&i_\textbf{s}
\beta^*_{\textbf{t},\textbf{s}}i^*_\textbf{t}i_\textbf{t}
u_\textbf{t}~~(\mbox{as}~~i_\textbf{t}\beta_{\textbf{t},\textbf{s}}=i_\textbf{s})\\&=&
i_\textbf{s}u_\textbf{s}. \Eea We first claim that the bounded net
$\{i_\textbf{s}u_\textbf{s}:\textbf{s}\in J_t\}$ converges to a
vector $v_t.$  For $\textbf{s}\leq\textbf{t}\in J_t,$ and $a\in
\mathcal E_t,$ \Bea |\langle
i_\textbf{t}u_\textbf{t}-i_\textbf{s}u_\textbf{s},a\rangle| &=&
|\langle(I_{\mathcal
E_t}-i_\textbf{s}i^*_\textbf{s})i_\textbf{t}u_\textbf{t},a\rangle|\\
&\leq& \|u_\textbf{t}\|\|(I_{\mathcal
E_t}-i_\textbf{s}i^*_\textbf{s})a\|\\ &\leq& e^{kt}\|(I_{\mathcal
E_t}-i_\textbf{s}i^*_\textbf{s})a\|. \Eea  As the net of projections
$\{i_\textbf{s}i^*_\textbf{s}: \textbf{s} \in J_t \}$ converges
strongly to identity, $(\langle
i_\textbf{s}u_\textbf{s},a\rangle)_{\textbf{s}\in J_t}$ is a Cauchy
net. Set $\phi (a)=\mbox{lim}_{\textbf{s}\in J_t}\langle
i_\textbf{s}u_\textbf{s},a\rangle.$ So there exists a unique vector
$v_t\in \mathcal E_t$ such that
$$\phi(a)=\mbox{lim}_{\textbf{s}\in J_t}\langle i_\textbf{s}u_\textbf{s},a\rangle=\langle v_t,a\rangle.$$
Also $$ \langle i_\textbf{s}i^*_\textbf{s}v_t,a\rangle = \langle
v_t, i_\textbf{s}i^*_\textbf{s}a\rangle =
\mbox{lim}_{\textbf{t}\in J_t}\langle
i_\textbf{t}u_\textbf{t},i_\textbf{s}i^*_\textbf{s}a\rangle=
\mbox{lim}_{\textbf{t}\in J_t}\langle
i_\textbf{s}\beta^*_{\textbf{t},\textbf{s}}u_\textbf{t},a\rangle=
\langle i_\textbf{s}u_\textbf{s},a\rangle . $$ So we get that for
any $\textbf{s}\in J_t,$
$$i_\textbf{s}i^*_\textbf{s}v_t=i_\textbf{s}u_\textbf{s}.$$  As
$\{i_\textbf{s}i^*_\textbf{s}:\textbf{s}\in J_t\}$ converges
strongly to identity of $\mathcal E_t,$ it proves that
$(i_\textbf{s}u_\textbf{s})_{\textbf{s}\in J_t}$ converges to $v_t$
in Hilbert space norm.

Now we claim that $(v_t)_{t>0}$ is a unit of the product system $(\mathcal E,B).$ For $a_1, a_2, \ldots a_k$ in
$E_s$ and $b_1, b_2, \ldots b_k$ in $E_t$ ($k\geq 1),$
 \Bea\langle B_{s,t}v_{s+t},\sum a_i\otimes b_i\rangle &=& \sum\langle v_{s+t},B^*_{s,t}(a_i\otimes b_i)\rangle \\&=& \sum
\lim_{\textbf{s}\smile \textbf{t}\in J_s\smile  J_t}\langle
i_{\textbf{s}\smile \textbf{t}}(u_\textbf{s}\otimes u_\textbf{t}),
B^*_{s,t}(a_i\otimes b_i)\rangle \\&=& \sum \lim _{\textbf{s}\smile \textbf{t}\in J_s\smile J_t} \langle
i_\textbf{s}u_\textbf{s}\otimes i_\textbf{t}u_\textbf{t},a_i\otimes
b_i\rangle\\&=& \sum \lim _{\textbf{s}\in J_s} \langle
i_\textbf{s}u_\textbf{s},a_i\rangle \lim_{\textbf{t}\in J_t}\langle
i_\textbf{t}u_\textbf{t},b_i\rangle\\ &=& \sum \langle
v_s,a_i\rangle\langle v_t,b_i\rangle\\ &=& \langle v_s\otimes
v_t,\sum a_i\otimes b_i\rangle . \Eea  This  proves that $v$ is a
unit. Finally, for $a\in E_\textbf{t}$,  we have \Bea \langle i^*_tv_t,a\rangle &=& \langle v_t,i_ta\rangle\\
&=& \lim _\textbf{r}\langle i_\textbf{r}u_\textbf{r},i_ta\rangle\\ &=&
\lim _\textbf{r}\langle i^*_ti_\textbf{r}u_\textbf{r},a\rangle\\ &=&
\lim _\textbf{r}
\langle\beta^*_{t,\textbf{r}}i^*_\textbf{r}i_\textbf{r}u_\textbf{r},a\rangle\\
&=&
\lim _\textbf{r}\langle\beta^*_{t,\textbf{r}}u_\textbf{r},a\rangle\\&=&\langle
u_t,a\rangle . \Eea which implies $i^*_tv_t=u_t.$ \qed

By this Theorem we see that if $u$ is a unit of an inclusion system
$(E, \beta )$ there exists a unique unit $\hat {u}$ in $(\mathcal E,
B)$ such that $i^*(\hat{u})=u$. We call $\hat {u}$ as the `lift' of
$u$. It is to be noted that for two units $u, v$ of the inclusion
system,
 $\langle \hat {u}_t, \hat{v}_t\rangle = \lim _{\textbf{s}\in J_t}\langle u_{\textbf{s}}, v_{\textbf{s}}\rangle.$ This helps us to compute
 covariance functions \cite{Ar01} of units.

\begin{thm}
Let $(E, \beta ), (F, \gamma)$ be two inclusion systems generating
two product systems $(\mathcal E, B), (\mathcal F, C)$ respectively.
Let $i, j$ be the respective inclusion maps. Suppose $A: (E, \beta
)\to (F, \gamma)$ is a  weak morphism then there exists a unique
morphism $\hat {A }: (\mathcal E, B)\to (\mathcal F, C)$ such that
$A_s= j_s^*\hat{A_s}i_s$ for all $s.$ This is a one to one
correspondence of weak morphisms. Further more, $\hat {A}$ is
isometric/unitary if $A$ is isometric/unitary.

\end{thm}

Proof:  If $\hat{A}$ is a morphism of product systems then $\{A_s=j^*_s\hat{A_s}i_s\}_{s>0}$ is clearly a weak
morpism of inclusion systems. Conversely  suppose $A: (E, \beta )\to (F, \gamma)$ is a morphism with
$\|A_t\|\leq e^{kt}$ for some $k>0.$ Define $A_\textbf{s}:E_\textbf{s}\rightarrow F_\textbf{s}$  by
$A_\textbf{s}=A_{s_1}\otimes \cdots \otimes A_{s_n}$. Let $i_\textbf{s}: E_\textbf{s}\rightarrow \mathcal E_s$
and $j_\textbf{s}:F_\textbf{s}\rightarrow \mathcal F_s$ be the canonical maps. The hypothesis implies that for
$\textbf{s} \leq \textbf{t}$ ,
$$\gamma^*_{\textbf{s},\textbf{t}}A_\textbf{t}\beta_{\textbf{s},\textbf{t}}=A_\textbf{s}.$$
Consider for $\textbf{s} \in J_s$,
$\Phi_\textbf{s}=j_\textbf{s}A_\textbf{s}i^*_\textbf{s}$. Set
$P_\textbf{r}=j_\textbf{r}j^*_\textbf{r}$ and
$Q_\textbf{r}=i_\textbf{r}i^*_\textbf{r}$. A simple computation
shows that for $\textbf{r}\leq \textbf{s}$,
$$P_\textbf{r}\Phi_\textbf{s}Q_\textbf{r}=\Phi_\textbf{r}.$$
For $\textbf{s}\leq\textbf{t}\in J_t,$ $a\in\mathcal E_t$ and
$b\in\mathcal F_t$ \Bea |\langle
(\Phi_\textbf{t}-\Phi_\textbf{s})a,b\rangle| &=& |\langle
(\Phi_\textbf{t}-P_\textbf{s}\Phi_\textbf{t}Q_\textbf{s})a,b\rangle|
\\ &\leq& |\langle
(\Phi_\textbf{t}-\Phi_\textbf{t}Q_\textbf{s})a,b\rangle| + |\langle
(\Phi_\textbf{t}Q_\textbf{s}-P_\textbf{s}\Phi_\textbf{t}Q_\textbf{s})a,b\rangle|
\\ &=& |\langle(I_{\mathcal
F_t}-Q_\textbf{s})a,\Phi^*_\textbf{t}b\rangle| + |\langle
\Phi_\textbf{t}Q_\textbf{s}a,(I_{\mathcal
E_t}-P_\textbf{s})b\rangle| \\& \leq & e^{kt}\|(I_{\mathcal
F_t}-Q_\textbf{s})a\|\|b\| + e^{kt}\|a\|\|b\|. \Eea  Imitating the
proof in the Theorem \ref{unit}, $(\Phi_\textbf{s})_{\textbf{s}\in
J_t}$ has a weak limit say $\hat{A_s}$. Now for $\textbf{s} \in
J_s$, we get \Bea \langle j^*_\textbf{s}\hat{A_s}i_\textbf{s}a ,
b\rangle &=& \langle \hat{A_s}i_\textbf{s}a,j_\textbf{s}b \rangle \\
&=& lim\langle \Phi_\textbf{r}i_\textbf{s}a,j_\textbf{s}b\rangle \\
&=& lim\langle
j^*_\textbf{s}j_\textbf{r}A_\textbf{r}i^*_\textbf{r}i_\textbf{t}a,b\rangle
\\ &=& lim\langle
\gamma^*_{\textbf{s},\textbf{r}}A_\textbf{r}\beta_{\textbf{s},\textbf{r}}a,b\rangle
\\ &=& \langle A_\textbf{s}a,b\rangle . \Eea  This implies that
$A_\textbf{s}= j_\textbf{s}^*\hat{A_s}i_\textbf{s}$ and in
particular $A_s= j_s^*\hat{A_s}i_s .$ Now we claim  that $\hat{A_s}$
is a morphism of product systems. For any $\textbf{s}\in J_s,$
$\textbf{t}\in J_t,$ $a\in E_\textbf{s},$ $b\in E_\textbf{t},$ $c\in
F_\textbf{s},$ $d\in F_\textbf{t}$ we have  \Bea \langle
C^*_{s,t}(\hat{A_s}\otimes\hat{A_t})B_{s,t}i_{\textbf{s}\smile \textbf{t}}(a\otimes
b),j_{\textbf{s}\smile \textbf{t}}(c\otimes d)\rangle &=& \langle
(\hat{A_s}\otimes\hat{A_t})(i_\textbf{s}\otimes
i_\textbf{t})(a\otimes b),j_\textbf{s}\otimes j_\textbf{t}(c\otimes
d)\rangle\\ &=& \langle j^*_\textbf{s}\hat{A_s}i_\textbf{s}\otimes
j^*_\textbf{t}\hat{A_t}i_\textbf{t}(a\otimes b), (c\otimes
d)\rangle\\ &=& \langle (A_\textbf{s}\otimes A_\textbf{t})(a\otimes
b), (c\otimes d)\rangle \\&=& \langle
A_{\textbf{s}\smile \textbf{t}}(a\otimes b), (c\otimes d)\rangle \\&=&
\langle
j^*_{\textbf{s}\smile \textbf{t}}\hat{A}_{s+t}i_{\textbf{s}\smile \textbf{t}}(a\otimes b), (c\otimes d)\rangle\\
&=& \langle \hat{A}_{s+t}i_{\textbf{s}\smile \textbf{t}}(a\otimes
b),j_{\textbf{s}\smile \textbf{t}}(c\otimes d)\rangle . \Eea The one to
one property can be proved imitating the proof in Theorem 10. The
second statement is obvious. \qed

As a special case we have the following universal property for
(weak/strong) morphisms.

\begin{cor}
Let $(\mathcal E, B)$ be a product system generated by an inclusion
system $(E, \beta)$ with canonical map $i$. Suppose $(\mathcal F,
C)$ is a product system with isometric morphisms of inclusion system
$m:E \rightarrow \mathcal F$. Then  there exists unique isometric
morphism of product system $\hat{m}:\mathcal E \rightarrow \mathcal
F$ such that $\hat {m}_si_s=m_s$ for all $s>0$.
\end{cor}

With basic theory of inclusion systems and their morphisms in place,
we look at inclusion systems arising from quantum dynamical
semigroups. Though this is part of folklore, as we are going to need
it in the next Section we put in some details. Let $H$ be a a
Hilbert space and let $\mathcal B(H)$ be the algebra of all bounded
operators on $H$. Let $\tau =\{ \tau _t: t\geq 0\}$ be a quantum
dynamical semigroup on $\mathcal B(H)$, that is, a one parameter
semigroup of normal, contractive, completely positive maps of
$\mathcal B(H).$ For $t\geq 0,$ let $(\pi _t, V_t, K_t)$ be a
 Stinespring dilation of $\tau _t$: $K_t$ is a
Hilbert space, $V_t\in \mathcal B(H, K_t),$ and $\pi _t$ is a normal
representation of $\mathcal B(H)$ on $K_t,$ such that,  $$ \tau
_t(X)=V_t^*\pi _t(X)V_t ~\forall X\in \mathcal B(H).$$
We will not need minimality ( $K_t= \overline
{\mbox {span}}\{\pi _t(X)V_th : X\in \mathcal B(H), h\in H\}).$ Now fix a
unit vector $a\in  H,$ take
$$E_t= \overline {\mbox {span}}\{ \pi _t(|a\rangle \langle g|)V_th: g,h
\in H\}\subseteq K_t.$$ Up to unitary equivalence the  Hilbert space $E_t$ does not depend upon the Stinespring
dilation or the choice of the reference vector $a.$ For any two unit vectors $a$ and $a^\prime,$ The map
$$\pi_t(|a\rangle\langle g|)V_th\rightarrow \pi_t(|a^\prime\rangle\langle g|)V_th $$ extends as a unitary between $E^a_t$ and
$E^{a^\prime}_t$ as
$$\big{\langle} \pi _t(|a\rangle \langle g_1|)V_th_1, \pi _t(|a\rangle \langle g_2|)V_th_2\big{\rangle}
= \big{\langle} h_1 , \tau _t(|g_1\rangle \langle g_2|)h_2\big{\rangle} = \big{\langle} \pi _t(|a^\prime\rangle
\langle g_1|)V_th_1, \pi _t(|a^\prime\rangle \langle g_2|)V_th_2\big{\rangle}$$ We may also construct $E_t$ more
abstractly by the usual quotienting and completing procedure on defining
$$ \big{\langle} g_1\otimes h_1 , g_2\otimes h_2\big{\rangle} _{\tau _t}:=
 \big{\langle} h_1 , \tau _t(|g_1\rangle \langle g_2|)h_2\big{\rangle} ,$$
 for $g_1\otimes h_1, g_2\otimes h_2 \in H^*\otimes H.$
Now fix an ortho-normal basis $\{
e_k\}$ of $H$ and define $\beta _{s, t}: E_{s+t}\to E_s\otimes E_t$
by
$$\beta _{s, t}(\pi _{s+t}(|a\rangle \langle g|)V_{s+t}h)
=\sum _k \pi _s(|a\rangle \langle g|)V_se_k\otimes \pi _t(|a\rangle \langle e_k|)V_th.$$ Then $(E, \beta )$ is
an inclusion system. Indeed,  \Bea & & \big{\langle} \sum _k \pi _s(|a\rangle \langle g_1|)V_se_k\otimes \pi
_t(|a\rangle \langle e_k|)V_th_1, \sum _l \pi _s(|a\rangle \langle g_2|)V_se_l\otimes \pi _t(|a\rangle
\langle e_l|)V_th_2 \big{\rangle}\\
&=& \sum _{k,l} \big{\langle} e_k,  \tau _s(|g_1\rangle\langle g_2|)e_l\big{\rangle} .
\big{\langle} h_1, \tau _t(|e_k\rangle \langle e_l|)h_2\big{\rangle} \\
&=& \big{\langle} h_1 , \tau _t ( |\sum _{k,l} \langle e_k,  \tau _s(|g_1\rangle\langle g_2|)e_l\rangle  e_k\rangle \langle e_l|)h_2\big{\rangle} \\
&=& \big{\langle} h_1 , \tau_t ( |\sum_l \tau_s(|g_1\rangle\langle g_2|)e_l\rangle\langle e_l|)h_2 \big{\rangle} \\&=& \big{\langle} h_1 , \tau_t(\tau_s(|g_1\rangle\langle g_2|)\sum_{l=1}^\infty(|e_l\rangle\langle e_l|))h_2 \big{\rangle} \\ &=& \big{\langle} h_1 , \tau _t (\tau _s(|g_1\rangle \langle g_2|))h_2\big{\rangle} \\
& =& \big{\langle} h_1 , \tau _{s+t}(|g_1\rangle \langle g_2|)h_2\big{\rangle} . \Eea So $\beta _{s,t}$ is an
isometry, and the associativity property can also be verified by direct computation: For $r, s, t>0,$
 \Bea
 & &
 (I_r\otimes \beta _{s+t})\beta _{r, s+t}(\pi _{r+s+t}(|a\rangle \langle
 g|)V_{r+s+t}g)\\
 &=& (I_r\otimes \beta _{s+t})\sum _k[\pi _{r}(|a\rangle \langle g|)V_re_k\otimes
 \pi _{s+t}(|a\rangle \langle e_k|)V_{s+t}g]\\
 & =& \sum _{k,l} \pi _r(|a\rangle \langle g|)V_re_k\otimes \pi _s(|a\rangle \langle e_k|)V_se_l
 \otimes \pi _t(a\rangle \langle e_l|)V_tg\\
 & =& (\beta _{r,s}\otimes I_t). (\pi _{r+s}(|a\rangle \langle g|)V_{r+s}e_l)\otimes \pi _t(|a\rangle \langle e_l|)V_tg\\
 &=& (\beta _{r,s}\otimes I_t)\beta _{r+s, t}(\pi _{r+s+t}(|a\rangle \langle g|)V_{r+s+t}g).
 \Eea
Now we will show that $\beta $ does not depend upon the choice of the orthonormal
 basis. Let $e=(e_i)_{i=1}^\infty $ and $f=(f_j)_{j=1}^\infty $ be two orthonormal bases of the Hilbert space $H.$ Now denoting the associated $\beta $ maps by $\beta ^e, \beta ^f$ respectively, we get
 \Bea & &  \big{\langle}\beta^e_{s,t}\pi_{s+t}(|a\rangle\langle g|)V_{s+t}h,  \beta^f_{s,t}\pi_{s+t}(|a\rangle\langle g|)V_{s+t}h \big{\rangle}\\ &=& \sum_{i,j}\big{\langle} \pi_s(|a\rangle\langle g|)V_se_i \otimes \pi_t(|a\rangle\langle e_i|)V_th , \pi_s(|a\rangle\langle g|)V_sf_j \otimes \pi_t(|a\rangle\langle f_j|)V_th \big{\rangle} \\ &=& \sum _{i,j} \big{\langle} e_i,  \tau _s(|g\rangle\langle g|)f_j\big{\rangle} .
\big{\langle} h, \tau _t(|e_i\rangle \langle f_j|)h\big{\rangle} \\ &=& \big{\langle} h , \tau _t ( |\sum
_{i,j} \langle e_i,  \tau _s(|g\rangle\langle g|)f_j\rangle  e_i\rangle \langle f_j|)h\big{\rangle} \\ &=&
\big{\langle} h , \tau_t(\tau_s(|g\rangle\langle g|)\sum_{j=1}^\infty(|f_j\rangle\langle f_j|))h \big{\rangle}
\\ &=& \big{\langle} h , \tau_{s+t}(|g\rangle\langle g|)h
\big{\rangle} . \Eea  So \Bea & &
\|\beta^e_{s,t}\pi_{s+t}(|a\rangle\langle g|)V_{s+t}h  -
\beta^f_{s,t}\pi_{s+t}(|a\rangle\langle g|)V_{s+t}h \|^2 \\ &=&
\big{\langle} h , \tau_{s+t}(|g\rangle\langle g|)h \big{\rangle} -
\big{\langle} h , \tau_{s+t}(|g\rangle\langle g|)h \big{\rangle} \\ & & -
\big{\langle} h , \tau_{s+t}(|g\rangle\langle g|)h \big{\rangle}  +
\big{\langle} h , \tau_{s+t}(|g\rangle\langle g|)h \big{\rangle} \\
&=& 0 \Eea Now we recall the dilation theorem for quantum dynamical
 semigroups (This was proved in \cite{Bh01} for unital quantum dynamical semigroups and was
 extended to the non-unital case in \cite{Bh02}):
Given a quantum dynamical semigroup $\tau $ on $\mathcal B(\mathcal
H)$ there exists a pair $(\theta , \mathcal K)$ where $\mathcal K$
is a Hilbert space containing $\mathcal H$ and $\theta $ is an
$E$-semigroup of $\mathcal B(\mathcal K)$ such that,
$$\tau _t(X)= P\theta _t(X)P ~~\forall ~X\in \mathcal B(\mathcal H), t\geq 0,$$
where $P$ is the projection of $\mathcal K$ onto $\mathcal H$ and
$X\in \mathcal B(\mathcal H)$ is identified with $PXP$ in $\mathcal
B(\mathcal K).$ Furthermore, we can choose $\mathcal K$ such that,
$$\overline {\mbox {span}}\{\theta _{r_1}(X_1)\ldots \theta
_{r_n}(X_n)h : r_1\geq r_2\geq \cdots \geq r_n\geq 0, X_1, \ldots
X_n\in \mathcal B(\mathcal H ), h\in \mathcal H, n\geq 0\}=\mathcal
K.$$ Such a pair $(\theta , \mathcal K)$  is unique up to unitary
equivalence and is called the minimal dilation of $\tau .$ The
minimal dilation $\theta $ is unital if and only if $\tau $ is
unital. In the following, we need another basic property of minimal
dilation: The vector $\theta _{r_1}(X_1)\ldots \theta _{r_n}(X_n)h$
appearing above remains unchanged if we drop any $\theta
_{r_k}(X_k)$ from the expression, if $X_k=1_{\mathcal H}$. In the
unital case, this fact follows easily  from the property that $P=
1_{\mathcal H}$ is an increasing projection for $\theta$. It is a
bit more involved in the non-unital case.

\begin{thm}
Let $\tau $ be a quantum dynamical semigroup of $\mathcal B(H)$ and let $(E, \beta )$ be the associated inclusion system defined above.
 Let $\theta$ acting
on $\mathcal B(K)$  (with $H\subset K$) be the minimal  E-semigroup
dilation of $\tau .$ Let $(\mathcal F, C)$ be the inclusion system
of $\theta $, considered as a quantum dynamical semigroup. Then
$(\mathcal F, C)$ is a product system and  is isomorphic to the
product system $(\mathcal E, B)$ generated by $(E, \beta )$.
\end{thm}

Proof: Since for every $t$, $\theta _t$ is a $*$-endomorphism, its
minimal Stinespring dilation is itself. Then it is easily seen that
$F_t= \overline {\mbox {span}}\{\theta _t (|a\rangle \langle x|)y :
x, y\in K\}$ and $C_{s, t}: {\mathcal F}_{s+t}\to {\mathcal
F}_s\otimes {\mathcal F}_t$, defined by
$$C_{s,t}(\theta _{s+t}(|a\rangle \langle x|)y= \sum _k\theta _s(|a\rangle
\langle x|)e_k\otimes \theta _t(|a\rangle \langle e_k|)y,$$
 is a unitary with
$$C_{s, t}^*(\theta _s(|a\rangle \langle x_1|)y_1\otimes \theta _t(|a\rangle \langle x_2|)y_2)
=\theta _{s+t}(|a\rangle \langle x_1|)\theta _t(|y_1\rangle \langle
x_2|)y_2.$$ Therefore $(\mathcal F, C)$ is a product system.

Define $m_s:E_s\rightarrow\mathcal F_s$ by $m_s(\pi _s(|a\rangle
\langle g|)V_sh)=\theta_s(|a\rangle \langle g|)h$ where $a$ is a
unit vector in $\mathcal H$. Clearly $m_s$ is a linear isometry.  We
see that $m$ is a strong morphism of inclusion systems: \Bea & &
C_{s,t}^*(m_s\otimes m_t)\beta _{s,t}(\pi _{s+t}(|a\rangle \langle
g|)V_{s+t}h)\\
& =&  C_{s,t}^*(m_s\otimes m_t)(\sum _k\pi _s(|a\rangle
\langle g|V_se_k\otimes \pi _t(|a\rangle \langle e_k|)V_th) \\
&=& C_{s,t}^*(\sum _k\theta _s(|a\rangle \langle g|)e_k\otimes
\theta _t(|a\rangle \langle e_k|)h)\\
&=& \sum _k\theta _{s+t}(|a\rangle \langle g|)\theta _t(|e_k\rangle
\langle e_k|)h\\
&=&\theta _{s+t}(|a\rangle \langle g|)\theta _t(1_{\mathcal H})h\\
&=&\theta _{s+t}(|a\rangle \langle g|)h\\
&=&m_{s+t}(\pi _{s+t}(|a\rangle \langle g|)V_{s+t}h).\Eea
 Now by Corollary 12, there exists an isometric morphism $\hat {m}:(\mathcal E,B
)\to (\mathcal F, C)$ satisfying $\hat {m}_si_s=m_s$ where $i: (E,
\beta )\to (\mathcal E, B)$ is the natural inclusion map. Unitarity
of $\hat {m}$ follows easily from the minimality of the dilation.
\qed


This theorem provides us with a plenty of inclusion systems with
finite dimensional systems, as we can consider contractive CP
semigroups on ${\mathcal B}(H)$ with dim $H<\infty .$  For instance,
one gets the inclusion system of  Example 2, by considering the following
CP semigroup on $\mathcal B (\mathbb C^2):$
$$T_t\left(
\begin{array}{cc}
a & b \\
c & d \\
\end{array}
\right)= e^{-\alpha t}\left(
\begin{array}{cc}
a+td & b \\
c & d\\
\end{array}\right)
,$$
where $\alpha $ is a suitable positive real number so as to make the semigroup contractive (This semigroup is not unital, but that does not matter). However, all
such inclusion systems would only generate type I product systems,
as  it is well-known that the associated product systems of CP
semigroups with bounded generators are always type I. This raises
the natural question as to whether inclusion systems $(E, \beta )$,
where dim $(E_t)\leq N$ for some natural number $N$ always generate
type I product systems  (One has to be a bit cautious here as the
product system generated may contain non-separable Hilbert spaces).
Recently this  has been answered in the affirmative by B. Tsirelson
for the case $N=2.$ (See \cite{Tsi01}, \cite{Tsi02}).

\begin{prop}
 Let $(\mathcal F, C)$ be a spatial product system and
let $(E, \beta )$ be the inclusion system formed by the linear spans
of units (See Example 2). Then the product system $(\mathcal E, B)$
generated by $(E, \beta )$ is the type I part of $(\mathcal F, C).$
\end{prop}

Proof: This is obvious, as the space $E_{\textbf{t}}=E_{t_n}\otimes
\cdots \otimes E_{t_1}$ can be identified with $\overline
{\mbox{span}}\{u^n_{t_n}\otimes \cdots \otimes u^1_{t_1}: u^1,
\ldots , u^n\in \mathcal U^{\mathcal E}\}.$

\qed

\section{amalgamation}

Suppose $H$ and $K$ are two Hilbert spaces and $D:K \to H$ is a linear contraction.
 Define a semi inner product on $H \oplus K$ by \Bea
\big{\langle} \left( \begin{array}{c} u_1 \\
                                      v_1 \\
                                    \end{array}
                                  \right),\left(
                                            \begin{array}{c}
                                              u_2 \\
                                              v_2 \\
                                            \end{array}
                                          \right)
\big{\rangle}_D &=& \big{\langle} u_1,u_2\big{\rangle} + \big{\langle} u_1,D
v_2 \big{\rangle} + \big{\langle} Dv_1,u_2 \big{\rangle} + \big{\langle} v_1,v_2 \big{\rangle} \\
& =& \big{\langle} \left(
\begin{array}{c}
u_1 \\
v_1 \\
\end{array}
\right) , \tilde {D}\left(
                                                                   \begin{array}{c}
                                                                     u_2 \\
                                                                     v_2 \\
                                                                   \end{array}
\right) \big{\rangle} , \\
\Eea
 where $\tilde {D} := \left[
                                                     \begin{array}{cc}
                                                       I & D \\
                                                       D^* & I \\
                                                     \end{array}
                                                   \right].$ Note that as $D$
is contractive, $\tilde {D}$  is positive definite.  Take
$$N=\{\left(
       \begin{array}{c}
         u \\
         v \\
       \end{array}
     \right)
:\big{\langle} \left(
          \begin{array}{c}
            u \\
            v \\
          \end{array}
        \right)
,\left(
   \begin{array}{c}
     u \\
     v \\
   \end{array}
 \right)
\big{\rangle}_D=0\}.$$ Then $N$ is the kernel of  bounded operator $\tilde{D}$ and hence it is a closed subspace of $H\oplus K$. Set
$G$ as completion of  $(H\oplus K)/N$ ) with respect to norm of $\big{\langle}. , .\big{\rangle} _D.$ We denote
$G$ by $H\oplus_D K$ and further denote the image of vector  $\left(
            \begin{array}{c}
              u \\
              v \\
            \end{array}
          \right)$
 by $\left[
            \begin{array}{c}
              u \\
              v \\
            \end{array}
          \right]$ for $u \in H$ and $v \in K .$
Now $$\big{\langle} \left[
                                                               \begin{array}{c}
                                                                 u_1 \\
                                                                 0 \\
                                                               \end{array}
                                                             \right]
          ,\left[
             \begin{array}{c}
               u_2 \\
               0 \\
             \end{array}
           \right]\big{\rangle}_D
          =\big{\langle} u_1,u_2 \big{\rangle}_H; ~~ \big{\langle} \left[
                                                               \begin{array}{c}
                                                                 0 \\
                                                                 v_1 \\
                                                               \end{array}
                                                             \right]
          ,\left[
             \begin{array}{c}
               0 \\
               v_2 \\
             \end{array}
           \right]\big{\rangle}_D
          =\big{\langle} v_1,v_2 \big{\rangle}_K. $$ So $H$ and $K$ are naturally embedded in $H\oplus_D K$ and their closed liner
span is $H\oplus_D K$ but they need not be orthogonal. We call
$H\oplus _DK$ as the amalgamation of $H$ and $K$ through $D$. It is
to be noted that if range~ $(\tilde {D})$ is closed, then no completion
is needed in the construction, and every vector of $G$ is of the
form $\left[
            \begin{array}{c}
              u \\
              v \\
            \end{array}
          \right]$ for $u \in H$ and $v \in K.$

In the converse direction, if $H$ and $K$ are two closed subspaces of a Hilbert space $G$. Then by a simple
application of Riesz representation theorem, there exists  unique contraction $D : K \rightarrow H$ such that
for $u\in H,$ $v\in K$ $$\big{\langle} u,v\big{\rangle}_G=\big{\langle} u , Dv \big{\rangle} .$$

Now we consider amalgamation at the level of inclusion systems. Let $(E , \beta )$ and $(F , \gamma )$ be two
inclusion systems. Let $D=\{ D_s: s>0\}$ be a weak contractive morphism from $F$ to $E$. Define
$G_s:=E_s\oplus_{D_s}F_s$ and $\delta_{s,t}:=i_{s,t}(\beta_{s,t}\oplus_D\gamma_{s,t})$ where
$i_{s,t}:(E_s\otimes E_t)\oplus_{D_s\otimes D_t}(F_s\otimes F_t)\rightarrow G_s\otimes G_t$ is the map defined
by
$$i_{s,t}\left[
                     \begin{array}{c}
                       u_1\otimes u_2 \\
                       v_1\otimes v_2 \\
                     \end{array}
                   \right]=\left[
                             \begin{array}{c}
                               u_1 \\
                               0 \\
                             \end{array}
                           \right]\otimes \left[
                                            \begin{array}{c}
                                              u_2 \\
                                              0 \\
                                            \end{array}
                                          \right] + \left[
                             \begin{array}{c}
                               0 \\
                               v_1 \\
                             \end{array}
                           \right]\otimes \left[
                                            \begin{array}{c}
                                              0 \\
                                              v_2 \\
                                            \end{array}
                                          \right],
$$
and $(\beta_{s,t}\oplus_D\gamma_{s,t}):E_{s+t}\oplus_{D_{s+t}} F_{s+t}\rightarrow E_s\otimes
E_t\oplus_{D_s\otimes D_t}F_s\otimes F_t$ is the map defined by $$(\beta_{s,t}\oplus_D\gamma_{s,t})\left[
                                                            \begin{array}{c}
                                                              u \\
                                                              v \\
                                                            \end{array}
                                                          \right]
= \left[
    \begin{array}{c}
      \beta_{s,t}(u) \\
      \gamma_{s,t}(v) \\
    \end{array}
  \right]
$$

\begin{lem}
The maps $i_{s,t}:(E_s\otimes E_t)\oplus_{D_s\otimes D_t}(F_s\otimes F_t)\rightarrow G_s\otimes G_t$ and
$(\beta_{s,t}\oplus_D\gamma_{s,t}):E_{s+t}\oplus_{D_{s+t}} F_{s+t}\rightarrow (E_s\otimes E_t)
\oplus_{D_s\otimes D_t} (F_s\otimes F_t)$
 are well defined isometries.
\end{lem}
Proof: It is enough to check that the maps are inner product preserving on elementary tensors. Observe that \Bea
& & \big{\langle} i_{s,t} \left[
            \begin{array}{c}
              u_1 \otimes u_2 \\
              v_1 \otimes v_2 \\
            \end{array}
          \right], i_{s,t} \left[
            \begin{array}{c}
              u^\prime_1 \otimes u^\prime_2 \\
              v^\prime_1 \otimes v^\prime_2 \\
            \end{array}
          \right] \big{\rangle}\\ &=& \big{\langle} \left[
                             \begin{array}{c}
                               u_1 \\
                               0 \\
                             \end{array}
                           \right]\otimes \left[
                                            \begin{array}{c}
                                              u_2 \\
                                              0 \\
                                            \end{array}
                                          \right] + \left[
                             \begin{array}{c}
                               0 \\
                               v_1 \\
                             \end{array}
                           \right]\otimes \left[
                                            \begin{array}{c}
                                              0 \\
                                              v_2 \\
                                            \end{array}
                                          \right] , \left[
                             \begin{array}{c}
                               u^\prime_1 \\
                               0 \\
                             \end{array}
                           \right]\otimes \left[
                                            \begin{array}{c}
                                              u^\prime_2 \\
                                              0 \\
                                            \end{array}
                                          \right] + \left[
                             \begin{array}{c}
                               0 \\
                               v^\prime_1 \\
                             \end{array}
                           \right]\otimes \left[
                                            \begin{array}{c}
                                              0 \\
                                              v^\prime_2 \\
                                            \end{array}
                                          \right] \big{\rangle}\\ &=& \big{\langle} u_1,u^\prime_1 \big{\rangle}\big{\langle} u_2,u_2^\prime\big{\rangle} +  \big{\langle} v_1,v_1^\prime \big{\rangle}\big{\langle} v_2,v^\prime_2\big{\rangle} + \big{\langle} u_1 , D_sv^\prime_1\big{\rangle} \big{\langle} u_2 ,
          D_tv^\prime_2\big{\rangle} + \big{\langle} D_sv_1 , u^\prime_1\big{\rangle} \big{\langle} D_tv_2 , u^\prime_2\big{\rangle} \\ &=& \big{\langle} \left[
            \begin{array}{c}
              u_1 \otimes u_2 \\
              v_1 \otimes v_2 \\
            \end{array}
          \right] ,  \left[
            \begin{array}{c}
              u^\prime_1 \otimes u^\prime_2 \\
              v^\prime_1 \otimes v^\prime_2 \\
            \end{array}
          \right] \big{\rangle}. \Eea
and
\Bea & &  \big{\langle}(\beta_{s,t}\oplus_D\gamma_{s,t}) \left[
                                                            \begin{array}{c}
                                                              u \\
                                                              v \\
                                                            \end{array}
                                                          \right] , (\beta_{s,t}\oplus_D\gamma_{s,t})  \left[
                                                            \begin{array}{c}
                                                              u^\prime \\
                                                              v^\prime \\
                                                            \end{array}
                                                          \right] \big{\rangle}\\ &=& \big{\langle} \left[
    \begin{array}{c}
      \beta_{s,t}(u) \\
      \gamma_{s,t}(v) \\
    \end{array}
  \right],\left[
    \begin{array}{c}
      \beta_{s,t}(u^\prime) \\
      \gamma_{s,t}(v^\prime) \\
    \end{array}
  \right]\big{\rangle} \\ &=&  \big{\langle} \beta_{s,t}(u),\beta_{s,t}(u^\prime)\big{\rangle} + \big{\langle} \gamma_{s,t}(v),\gamma_{s,t}(v^\prime) \big{\rangle}\\
   & & + \big{\langle} \beta_{s,t}(u),(D_s\otimes D_t)\gamma_{s,t}(v^\prime )\big{\rangle} + \big{\langle} (D_s\otimes D_t)\gamma_{s,t}
(v),\beta_{s,t}(u^\prime ) \big{\rangle} \\ &=& \big{\langle} u,u^\prime\big{\rangle} + \big{\langle} v,v^\prime
\big{\rangle}+\big{\langle} u,D_{s+t}v^\prime \big{\rangle} + \big{\langle} D_{s+t}v,u^\prime \big{\rangle} \\ &=&
\big{\langle} \left[
                             \begin{array}{c}
                               u \\
                               v \\
                             \end{array}
                           \right],\left[
                             \begin{array}{c}
                               u^\prime \\
                               v^\prime \\
                             \end{array}
                           \right]\big{\rangle} . \Eea \qed
\begin{prop}
Let $(G,\delta )= \{ G_s, \delta _{s,t}: s,t>0\}$ be defined as
above. Then \{$G,\delta$\} forms an inclusion system.
\end{prop}Proof:  Being composition of two isometries, $\delta _{s,t}$ is an
isometry. Define $i_{r,s,t}: (E_r\otimes $ $ E_s \otimes E_t) \oplus_{D_r\otimes D_s\otimes D_t} (F_r\otimes
F_s\otimes F_t) \rightarrow  G_r\otimes G_s\otimes G_t$ by $$i_{r,s,t}\left[
                                                         \begin{array}{c}
                                                           u_1\otimes u_2\otimes u_3 \\
                                                           v_1\otimes v_2\otimes v_3 \\
                                                         \end{array}
                                                       \right]=\left[
                                                                 \begin{array}{c}
                                                                   u_1 \\
                                                                   0 \\
                                                                 \end{array}
                                                               \right]\otimes \left[
                                                                                \begin{array}{c}
                                                                                  u_2 \\
                                                                                  0 \\
                                                                                \end{array}
                                                                              \right]\otimes \left[
                                                                                               \begin{array}{c}
                                                                                                 u_3 \\
                                                                                                 0 \\
                                                                                               \end{array}
                                                                                             \right]+ \left[
                                                                 \begin{array}{c}
                                                                   0 \\
                                                                   v_1 \\
                                                                 \end{array}
                                                               \right]\otimes \left[
                                                                                \begin{array}{c}
                                                                                  0 \\
                                                                                  v_2 \\
                                                                                \end{array}
                                                                              \right]\otimes \left[
                                                                                               \begin{array}{c}
                                                                                                 0 \\
                                                                                                 v_3 \\
                                                                                               \end{array}
                                                                                             \right]
$$

It can be shown similarly that $i_{r,s,t}$ is an isometry. For
$r,s,t \in \mb R_+$
\Bea (\delta_{r,s}\otimes 1_{G_t})\delta_{r+s,t} &=& [i_{r,s}(\beta_{r,s}\oplus_D\gamma_{r,s})\otimes
1_{G_t}]i_{r+s,t}(\beta_{r+s,t}\oplus_D\gamma_{r+s,t})\\ & = &  i_{r,s,t}((\beta_{r,s}\otimes
1_{E_t})\beta_{r+s,t}\oplus_D(\gamma_{r,s}\otimes 1_{E_t})\gamma_{r+s,t})  \Eea
Similarly
 $$(1_{G_r}\otimes\delta_{s,t})\delta_{r,s+t}=i_{r,s,t}((1_{E_r}\otimes\beta_{s,t})\beta_{r,s+t}\oplus_D(1_{F_r}\otimes\gamma_{s,t})\gamma_{r,s+t})$$

As $(E, \beta ), $ $(F ,\gamma )$ are inclusion systems, so is $(G ,
\delta ).$ \qed


\begin{defn}
The inclusion system $(G, \delta)$ constructed above is called the
amalgamation of inclusion systems $(E,\beta )$ and $(F, \gamma)$ via
the morphism $D$. If $(\mathcal E, B), (\mathcal F, C), $ and
$(\mathcal G, L)$ are product systems generated respectively by $(E,
\beta ), (F, \gamma ) ,$ and $(G, \delta )$, then $(\mathcal G, L)$
is said to be the amalgamated product of $(\mathcal E, B)$ and
$(\mathcal F, C)$ via $D$ and is denoted by $\mathcal {G} =:
\mathcal E\otimes _D\mathcal F .$

\end{defn}

In this Definition notice that as $D$ is a weak morphism of
inclusion system we will get a lift $\hat{D}: \mathcal F\to \mathcal
E .$ We can also define $\mathcal E\otimes_{\hat{D}}\mathcal F$. It
can be seen easily that product system generated by the amalgamated
product of $(E,\beta)$ and $(F,\gamma)$ via $D$ is same as $\mathcal
E\otimes_{\hat{D}}\mathcal F,$ so that the definition of amalgamated
product is unambiguous. This is true because of  the following
universal property of amalgamation.

\begin{prop} Let $(G, \delta )$ be the amalgamated inclusion system
of two inclusion systems $(E, \beta ), (F, \gamma)$ via a morphism
$D$ from $F$ to $E$. Let $(\mathcal G, L)$ be the product system
generated by $(G,\delta)$. Given any inclusion system $(H,\eta )$
with weak isometric morphisms $i:(E,\beta)\rightarrow (H, \eta )$
and $j:(F,\gamma )\rightarrow (H, \eta )$ with $\langle
i_ta,j_tb\rangle=\langle a,D_tb\rangle $ for all $a \in E_t$ and $b
\in F_t$, there exists unique isometric morphism of product system
$\hat{A}:(\mathcal G,L)\rightarrow (\mathcal H,W)$ such that
$l_t^*\hat{A}_t(k_t\left[
\begin{array}{c}
a \\
b \\ \end{array} \right] ) =(i_t(a)+j_t(b)),$ where $(\mathcal H,
W)$ is the product system generated by $(H, \eta ),$ and
$k_t:G_t\rightarrow \mathcal G_t,$ $l_t:H_t\rightarrow \mathcal H_t$
are respective  canonical maps.
\end{prop}
Proof: Define $A_t:G_t\rightarrow\mathcal H_t$ by $A_t\left[
                                                        \begin{array}{c}
                                                          a \\
                                                          b \\
                                                        \end{array}
                                                      \right]
=i_t(a)+j_t(b)$. Clearly $A_t$ is an isometry for each $t$. Now
$W^*_{s,t}(A_s\otimes A_t)\delta_{s,t}\left[
                                              \begin{array}{c}
                                                a \\
                                                b \\
                                              \end{array}
                                            \right]
=W^*_{s,t}(A_s\otimes A_t)i_{s,t}\left[
                                     \begin{array}{c}
                                       \beta_{s,t}(a) \\
                                       \gamma_{s,t}(b) \\
                                     \end{array}
                                   \right]
 =W^*_{s,t}(i_s\otimes
i_t)\beta_{s,t}(a)+W^*_{s,t}(j_s\otimes
j_t)\gamma_{s,t}(b)=i_{s+t}(a)+j_{s+t}(b)=A_{s+t}\left[
                                              \begin{array}{c}
                                                a \\
                                                b \\
                                              \end{array}
                                            \right]$. This means that $A$ is a weak
isometric morphism of inclusion systems. So it lifts to an isometric
morphism of product systems $\hat{A}:\mathcal G\rightarrow\mathcal
H_t$ such that $l_t^*\hat{A}_tk_t\left[
                                                        \begin{array}{c}
                                                          a \\
                                                          b \\
                                                        \end{array}
                                                      \right]=(i_t(a)+j_t(b))$ and the
proof is complete.\qed


Suppose $\phi =\{\phi _t : t\geq 0\}$ , $\psi =\{\psi _t: t\geq 0\}$
are CP semigroups on $\mathcal B(H)$ and $\mathcal B(K)$
respectively. Also suppose $\eta=\{\eta_t: t\geq 0\}$ is a family of
bounded operators on $B(K,H)$ such that $\tau=\{(\tau_t): t\geq 0\}$
defined by $$ \tau_t( \left(
\begin{array}{cc}
                       X & Y \\
                       Z & W
                     \end{array}
                   \right) = \left(
                              \begin{array}{cc}
                                \phi _t(X) & \eta_t(Y) \\
                                \eta_t(Z^*)^* & \psi _t(W)
                              \end{array} \right) $$  is a CP semigroup on $\mathcal B(H\oplus K)$. In particular,
this means that $\eta=\{ \eta_t :t\geq 0\}$ is a semigroup of
bounded maps  on $\mathcal B(K,H)$.

Let $(\pi_t, V_t, G_t)$ be the minimal Stinespring dilation of
$\tau_t$. Then clearly restrictions of $\pi_t$ to $\B(H)$ and $B(K)$
are dilations of $\phi _t$ and $\psi _t$. Suppose
$\hat{H_t}$=$\overline{\mbox{span}}\{\pi_t(X)V_th: X \in \B(H), h
\in H\}$ and $\hat{K_t}$=$\overline{\mbox{span}}\{\pi_t(Y)V_tk: Y
\in \B(K), k \in K\}$. Then $(\pi_t,V_t,\hat{H_t})$ and
$(\pi_t,V_t,\hat{K_t})$ are the minimal dilations of $\phi _t$ and
$\psi _t$ respectively. Fix unit vectors $a
\in H$ and $b \in K$ and $\left(%
\begin{array}{c}
  \frac{a}{\sqrt{2}} \\
  \frac{b}{\sqrt{2}} \\
\end{array}%
\right) \in H \oplus K $. Also fix orthonormal bases $\{e_i\}_{i\in I}$ of $H$ , $\{f_j\}_{j\in J}$ of $K,$
where $I,J$ are some indexing sets.  Then of course  $\{e_i\}_{i\in I}\cup \{f_j\}_{j\in J}$  is an orthonormal
basis of $H \oplus K$. Define inclusion systems $(E, \beta ),$  $(F,\gamma ), $ and  $(G,\delta )$ corresponding
to $\phi ,$  $\psi ,$  and $\tau .$ Define $D_t:F_t \rightarrow E_t$ by $D_t={P_{E_t}\pi_t(|a\rangle\langle
b|)|_{F_t}}$ (where $P_{E_t}$ is the projection onto $E_t ).$
\begin{thm}\label{amalgamation}
Let $\phi $ , $\psi $ , $\tau$ be CP semigroups and $(E,\beta)$ , $(F,\gamma)$ , $(G,\delta)$ be their
corresponding inclusion systems as above. Then $D=\{ D_t: t>0\}$  is a contractive morphism from $(F, \gamma )$
to $(E, \beta ).$ Moreover, $(G,\delta)$ is isomorphic to amalgamated sum of $(E, \beta)$ and $(F, \gamma)$ via
$D.$ \end{thm} Proof: Clearly each $D_t$ is contractive. To see that they form a morphism, we make the
following computations:\Bea & & \big{\langle} \beta_{s,t}\pi_{s+t}(|a\rangle\langle g|)V_{s+t}h,(D_s\otimes
D_t)\gamma_{s+t}\pi_{s+t}(|b\rangle \langle
g^\prime|)V_{s+t}h^\prime\big{\rangle} \\
&=&\sum_{i,j}\big{\langle} \pi_s(|a \rangle\langle g|)V_se_i\otimes \pi_t(|a\rangle\langle e_i)V_th,(D_s\otimes
D_t)\pi_s(|b\rangle\langle
g^\prime|)V_sf_j\otimes \pi_t(|b\rangle\langle f_j|)V_th^\prime\big{\rangle} \\
&=&\sum_{i,j}\big{\langle} \pi_s(|a\rangle\langle g|)V_se_i,P_{E_s}\pi_s(|a\rangle\langle
b|)\pi_s(|b\rangle\langle g^\prime|)V_sf_j \big{\rangle} \\
& & .\big{\langle} \pi_t(|a\rangle\langle
e_i|)V_th,P_{E_t}\pi_t(|a\rangle\langle b|)\pi_t(|b\rangle\langle f_j|)V_th^\prime \big{\rangle} \\
&=&\sum_{i,j}\big{\langle} e_i,\eta_s(|g\rangle\langle g^\prime|)f_j\big{\rangle} \big{\langle} h
,\eta_t(|e_i\rangle\langle f_j|)h^\prime \big{\rangle} \\
&=& \big{\langle} h , \eta _t( [\sum _{i,j}\big{\langle} e_i, \eta _s(|g\rangle \langle g^\prime )|f_j
\big{\rangle} |e_i\rangle \langle f_j|])h^\prime
\big{\rangle} \\
&=&\big{\langle} h,\eta_t(\eta_s(|g\rangle\langle
g^\prime|))h^\prime \big{\rangle} \\
&=&\big{\langle} h,\eta_{s+t}(|g\rangle\langle g^\prime|)h^\prime \big{\rangle} \\
&=&\big{\langle} \pi_{s+t}(|a\rangle\langle g|)V_{s+t}h,P_{E_{s+t}}\pi_{s+t}(|a\rangle\langle
b|)\pi_{s+t}(|b\rangle\langle g^\prime|)h^\prime \big{\rangle} . \Eea This proves the first part. Now define
$U_s:G_s\rightarrow E_s\oplus _{D_s}F_s$ by
$$U_s\pi_s(|\left(
\begin{array}{c}
 \frac{a}{\sqrt{2}}  \\
 \frac{b}{\sqrt{2}} \\
\end{array}
\right)\rangle\langle\left(
\begin{array}{c}
  g \\
  g^\prime \\
\end{array}
\right)|)V_s\left(
\begin{array}{c}
  h \\
  h^\prime \\
\end{array}
\right)=\left[
\begin{array}{c}
  \pi_s(|a\rangle\langle g|)V_sh \\
  \pi_s(|b\rangle\langle g^\prime|)V_sh^\prime \\
\end{array}
\right].$$ Clearly $U_s$ is linear and onto.
 Now \Bea   \parallel \left[
\begin{array}{c}
  \pi_s(|a\rangle\langle g|)V_sh \\
  \pi_s(|b\rangle\langle g^\prime|)V_sh^\prime \\
\end{array}
\right]\parallel^2
 &=& \big{\langle} h , \phi _s(|g\rangle\langle
g|)h\big{\rangle} + \big{\langle} h^\prime , \psi _s(|g^\prime\rangle\langle g^\prime|)h^\prime\big{\rangle} \\
& &  +
\big{\langle} h , \eta_s(|g\rangle\langle g^\prime|)h^\prime \big{\rangle} + \big{\langle}
\eta_s(|g\rangle\langle g^\prime|)h^\prime , h \big{\rangle} \\ &=& \big{\langle} \left(
     \begin{array}{c}
       h \\
       h^\prime \\
     \end{array}
   \right),
\left(
  \begin{array}{cc}
    \phi_s(|g\rangle\langle g^\prime|) & \eta_s(|g\rangle\langle g^\prime |) \\
    \eta_s(|g^\prime\rangle\langle g|^*)^* & \psi _s(|g^\prime\rangle\langle g^\prime|) \\
  \end{array}
\right) \left(
  \begin{array}{c}
    h \\
    h^\prime \\
  \end{array}
\right) \big{\rangle} \\
&=& \|\pi_s(|\left(%
\begin{array}{c}
  a/\sqrt{2} \\
  b/\sqrt{2} \\
\end{array}%
\right)\rangle \langle \left(%
\begin{array}{c}
  g \\
  g^\prime \\
\end{array}%
\right)|)V_s\left(%
\begin{array}{c}
  h \\
  h^\prime \\
\end{array}%
\right)\|^2 , \Eea
implies that $U_s$ is a unitary operator. In a
similar way, strong morphism property follows from:

$$ (U_s\otimes U_t)\delta_{s,t}\pi_{s+t}(|\left(%
\begin{array}{c}
  \frac{a}{\sqrt{2}} \\\frac{b}{\sqrt{2}} \\
\end{array}%
\right)\rangle \langle \left(%
\begin{array}{c}
  g \\
  g^\prime \\
\end{array}%
\right)|)V_{s+t}\left(%
\begin{array}{c}
  h \\
  h^\prime \\
\end{array}%
\right) $$
 \Bea
 &=& \sum_i
U_s\pi_s(|\left(%
\begin{array}{c}
  \frac{a}{\sqrt{2}} \\
  \frac{b}{\sqrt{2}} \\
\end{array}%
\right)\rangle \langle \left(%
\begin{array}{c}
  g \\
  g^\prime \\
\end{array}%
\right)|)V_s\left(%
\begin{array}{c}
 e_i\\
 0 \\
\end{array}%
\right) \otimes U_t \pi_t(|\left(%
\begin{array}{c} \frac{a}{\sqrt{2}} \\
  \frac{b}{\sqrt{2}} \\
\end{array}
\right)\rangle\langle\left(%
\begin{array}{c}
 e_i\\
 0 \\
\end{array}%
\right)|)V_t\left(%
\begin{array}{c}
  h\\
  h^\prime\\
\end{array}%
\right)\\
& &    +   \sum_j
U_s\pi_s(|\left(%
\begin{array}{c}
   \frac{a}{\sqrt{2}}\\
   \frac{b}{\sqrt{2}}\\
\end{array}
\right)\rangle\langle\left(%
\begin{array}{c}
  g \\
  g^\prime \\
\end{array}%
\right)|)V_s\left(%
\begin{array}{c}
 0\\
 f_j \\
\end{array}%
\right) \otimes  U_t\pi_t(|\left(%
\begin{array}{c}
  \frac{a}{\sqrt{2}} \\
  \frac{b}{\sqrt{2}} \\
\end{array}%
\right)\rangle\langle\left(%
\begin{array}{c}
 0\\
 f_j \\
\end{array}%
\right)|)V_t\left(%
\begin{array}{c}
  h\\
  h^\prime\\
\end{array}%
\right) \Eea
\Bea
 &=& \sum_i \left[
                    \begin{array}{c}
                      \pi_s(|a\rangle\langle g|)V_se_i \\
                      0 \\
                    \end{array}
                  \right]\otimes \left[
                                   \begin{array}{c}
                                     \pi_t(|a\rangle\langle e_i|)V_th \\
                                     0 \\
                                   \end{array}
                                 \right] \\
                                 & & + \sum_j \left[
                    \begin{array}{c}
                      0 \\
                      \pi_t(|b\rangle\langle g^\prime|)V_tf_j \\
                    \end{array}
                  \right]\otimes \left[
                                   \begin{array}{c}
                                     0 \\
                                     \pi_t(|b\rangle\langle f_j|)V_th^\prime \\
                                   \end{array}
                                 \right] \\ &=& i_{s,t}\left[
                                                             \begin{array}{c}
                                                               \sum_i\pi_s(|a\rangle\langle g|)V_se_i \otimes \pi_t(|a\rangle\langle e_i|)V_th \\
                                                               \sum_j\pi_s(|b\rangle\langle g^\prime|)V_sf_j \otimes \pi_t(|a\rangle\langle f_j|)V_th^\prime \\
                                                             \end{array}
                                                           \right] \\
&=& i_{s,t}\left[
                                              \begin{array}{c}
                                                \beta _{s,t}\pi_{s+t}(|a\rangle\langle g|)V_{s+t}h \\
                                                \gamma _{s,t}\pi_{s+t}(|b\rangle\langle g^\prime|)V_{s+t}h^\prime \\
                                              \end{array}
                                            \right] \\
& =& i_{s,t}(\beta _{s,t}\oplus _D \gamma _{s,t}) \left[
                                              \begin{array}{c}
                                                \pi_{s+t}(|a\rangle\langle g|)V_{s+t}h \\
                                                \pi_{s+t}(|b\rangle\langle g^\prime|)V_{s+t}h^\prime \\
                                              \end{array}
                                            \right] \\
&=& i_{s,t}(\beta_{s,t}\oplus _D\gamma_{s,t})U_{s+t}\pi_{s+t}(|\left(%
\begin{array}{c}
\frac{a}{\sqrt{2}} \\
  \frac{b}{\sqrt{2}} \\
\end{array}%
\right)\rangle\langle\left(%
\begin{array}{c}
  g \\
  g^\prime \\
\end{array}%
\right)|)V_{s+t}\left(%
\begin{array}{c}
  h \\
  h^\prime \\
\end{array}%
\right).  \Eea \qed

 Now we look at units of amalgamated products of
inclusion systems. \begin{lem}\label{formula}
 Let $(G, \gamma)$ be the amalgamated
product of  two inclusion systems $(E,\beta)$ and $(F,\gamma)$ via
$D$. Assume that range $(\tilde {D_t})$ is closed for every $t>0.$
Suppose $\left[
            \begin{array}{c}
              u_t \\
              v_t \\
            \end{array}
          \right]_{t>0}$  is a unit of $(G,\delta)$. Then $(u_s+D_s v_s)_{s>0}$ and
$(D^*_su_s+ v_s)_{s>0}$ are  units, provided they are non-trivial,  in $(E,\beta)$ and $(F,\gamma)$ respectively.

\end{lem}
Proof: As $\tilde{D_t}$ is closed ,  every vector of $G_t$ is given by $\left[
                                                                        \begin{array}{c}
                                                                          a \\
                                                                          b \\
                                                                        \end{array}
                                                                      \right]
$ for some $a\in E_t$ and $b\in F_t.$ Let $ \left[
               \begin{array}{c}
                 u_t \\
                 v_t \\
               \end{array}
             \right]
_{t>0} $ be a unit of the inclusion system $(G_t=E_t \oplus _D F_t, \delta_{s,t})$. So $$\left[
   \begin{array}{c}
     u_{s+t} \\
     v_{s+t} \\
   \end{array}
 \right]
=\delta^*_{s,t}\left[
                 \begin{array}{c}
                   u_s \\
                   v_s \\
                 \end{array}
               \right]
\otimes \left[
          \begin{array}{c}
            u_t \\
            v_t \\
          \end{array}
        \right]
=(\beta^*_{s,t}\oplus _D\gamma^*_{s,t})i^*_{s,t}\left[
                 \begin{array}{c}
                   u_s \\
                   v_s \\
                 \end{array}
               \right]
\otimes \left[
          \begin{array}{c}
            u_t \\
            v_t \\
          \end{array}
        \right].$$  Suppose $i^*_{s,t}\left[
                 \begin{array}{c}
                   u_s \\
                   v_s \\
                 \end{array}
               \right]
\otimes \left[
          \begin{array}{c}
            u_t \\
            v_t \\
          \end{array}
        \right]=\left[
                  \begin{array}{c}
                    z \\
                    w \\
                  \end{array}
                \right]
$, then we claim that $$  \left[ \begin{array}{cc} I&
D_s\otimes D_t \\
D_s^*\otimes D_t^* & I\end{array}\right] \left(
                  \begin{array}{c}
                    z \\
                    w \\
                  \end{array}
                \right) = \left(
                  \begin{array}{c}
                    (u_s+D_sv_s)\otimes (u_t+D_tv_t)\\
                    (D_s^*u_s+v_s)\otimes (D_t^*u_t\otimes v_t) \\
                  \end{array}
                \right). $$
This follows, as for $a\otimes a'\in E_s\otimes E_t$ and $b\otimes
b'\in F_s\otimes F_t,$ \Bea & & \langle
 \left[ \begin{array}{cc} I&
D_s\otimes D_t \\
D_s^*\otimes D_t^* & I\end{array}\right] \left(
                  \begin{array}{c}
                    z \\
                    w \\
                  \end{array}
                \right),
\left(
                  \begin{array}{c}
                    a\otimes a' \\
                    b\otimes b' \\
                  \end{array}
                \right) \rangle \\
&=& \langle  \left[
                  \begin{array}{c}
                    z \\
                    w \\
                  \end{array}
                \right],
\left[
                  \begin{array}{c}
                    a\otimes a'\\
                    b\otimes b'\\
                  \end{array}
                \right]\rangle \\
&=& \langle i_{s,t}^*( \left[
                  \begin{array}{c}
                    u_s \\
                    v_s \\
                  \end{array}
                \right]\otimes
                \left[
                  \begin{array}{c}
                    u_t\\
                    v_t \\
                  \end{array}
                \right] ),
\left[ \begin{array}{c} a\otimes a' \\ b\otimes b' \\ \end{array}
\right]\rangle \\
&=& \langle  \left[
                  \begin{array}{c}
                    u_s \\
                    v_s \\
                  \end{array}
                \right]\otimes
                \left[
                  \begin{array}{c}
                    u_t\\
                    v_t \\
                  \end{array}
                \right],
i_{s,t}\left[ \begin{array}{c} a\otimes a' \\ b\otimes b' \\
\end{array}
\right]\rangle \\
&=& \langle  \left[
                  \begin{array}{c}
                    u_s \\
                    v_s \\
                  \end{array}
                \right]\otimes
                \left[
                  \begin{array}{c}
                    u_t\\
                    v_t \\
                  \end{array}
                \right] ,
\left( \left[ \begin{array}{c} a \\ 0 \\ \end{array} \right] \otimes
\left[ \begin{array}{c} a' \\ 0 \\ \end{array} \right] + \left[
\begin{array}{c} 0 \\ b \\ \end{array} \right]
\otimes \left[ \begin{array}{c} 0 \\ b' \\ \end{array}
\right]\right) \rangle \\
&=& \langle \left(
                  \begin{array}{c}
                    (u_s+D_sv_s)\otimes (u_t+D_tv_t)\\
                    (D_s^*u_s+v_s)\otimes (D_t^*u_t\otimes v_t) \\
                  \end{array}
                \right),
\left(
                  \begin{array}{c}
                    a\otimes a'\\
                    b\otimes b' \\
                  \end{array}
                \right)\rangle .\Eea
Now for $c\in E_{s+t}, d\in F_{s+t}$ \Bea \langle\left(
\begin{array}{c}
u_{s+t}+D_{s+t}v_{s+t} \\
D_{s+t}^*u_{s+t}+v_{s+t} \\
\end{array}
\right),\left(
\begin{array}{c}
c \\
d \\
\end{array}
\right)\rangle &=&\langle\left[
\begin{array}{c}
u_{s+t} \\
v_{s+t} \\
\end{array}
\right],\left[
\begin{array}{c}
c \\
d \\
\end{array}
\right]\rangle \\
&=& \langle (\beta _{s,t}\oplus _D \gamma _{s,t})^*
i_{s,t}^* \left( \left[ \begin{array}{c} u_s\\ v_s \\ \end{array}
\right] \otimes \left[ \begin{array}{c} u_t \\ v_t \\ \end{array}
\right] \right), \left[ \begin{array}{c} c \\ d \\ \end{array}
\right]
\rangle \\
&=& \langle \left[ \begin{array}{c} z \\ w \\ \end{array} \right] ,
\left[ \begin{array}{c} \beta _{s,t}(c) \\ \gamma _{s,t}(d) \\
\end{array} \right] \rangle \\
&=& \langle
 \left[ \begin{array}{cc} I&
D_s\otimes D_t \\
D_s^*\otimes D_t^* & I\end{array}\right] \left(
                  \begin{array}{c}
                    z \\
                    w \\
                  \end{array}
                \right),
\left(
                  \begin{array}{c}
                    \beta_{s,t}(c) \\
                    \gamma _{s,t}(d) \\
                  \end{array}
                \right) \rangle \\
&=& \langle \left(
                  \begin{array}{c}
                    (u_s+D_sv_s)\otimes (u_t+D_tv_t)\\
                    (D_s^*u_s+v_s)\otimes (D_t^*u_t\otimes v_t) \\
                  \end{array}
                \right),
\left(
                  \begin{array}{c}
                    \beta _{s,t}(c)\\
                    \gamma _{s,t} (d) \\
                  \end{array}
                \right)\rangle \\
&=& \langle \left(
                  \begin{array}{c}
                    \beta _{s,t}^*((u_s+D_sv_s)\otimes (u_t+D_tv_t))\\
                    \gamma _{s,t}^*((D_s^*u_s+v_s)\otimes (D_t^*u_t\otimes v_t)) \\
                  \end{array}
                \right),
\left(
                  \begin{array}{c}
                    c\\
                    d \\
                  \end{array}
                \right)\rangle \Eea

Further  $\{(u_s+D_s v_s): s>0\}$ and $\{(D^*_su_s+ v_s): s>0\}$ are
exponentially bounded as $\parallel(u_s+D_s v_s)\parallel \leq
\parallel \left[
                               \begin{array}{c}
                                 u_s \\
                                 v_s \\
                               \end{array}
                             \right]
\parallel_{G_s}$ and  similarly $\parallel(D^*_su_s+ v_s)\parallel \leq
\parallel \left[
                               \begin{array}{c}
                                 u_s \\
                                 v_s \\
                               \end{array}
                             \right]\parallel_{G_s}$. Therefore they
                             are units in the corresponding inclusion systems. \qed


The Theorem \ref{amalgamation}{~} answers Powers question, even for general corners (not just those given by units), by showing that the product system of the CP semigroup formed is the amalgamation. Lemma \ref{formula} helps us in computing the units of such amalgamations. We also would like to compute Arveson index of the product systems. However, we have technical problem here. The index of a product system as defined by Arveson (\cite{Ar01}) needs measurability of units and we have not imposed any measurability structure on our inclusion systems or product systems. We do not intend to develop a theory of measurable inclusion systems here. Instead, we restrict ourselves to considering subsystems of measurable product systems, in the sense of Arveson.

Here after product systems we consider are Arveson systems and the units we consider are measurable.

Let  $\mathcal
U^\mathcal E$ denote the units of a product system  $\mathcal E$. Then the measurability ensures the existence a function
$$\gamma :\mathcal U^{\mathcal E}\times \mathcal U^{\mathcal E}\to \mathbb C$$
called the covariance function satisfying:
$$\langle u_t, v_t\rangle = e^{t\gamma (u, v) }~~~~\forall t,$$
for units $u, v.$ The function $\gamma $ is a conditionally positive definite function \cite{Ar01}. If $Z$ is a non-empty subset of ${\mathcal U}^{\mathcal E}$,  we may do the usual GNS construction for the kernel $\gamma $ restricted  to $Z\times Z$ to obtain a Hilbert space $H_Z$, which we call as the Arveson Hilbert space associated to $Z$. Note that the index of the product system $\mathcal E$ is nothing but the dimension of ${\mathcal K}:= H_{{\mathcal U}^{{\mathcal E}}}$ (Arveson Hilbert space of ${\mathcal U}^{\mathcal E})$. In \cite{Ar01}, it is shown that    there exists a bijection $u\mapsto (\lambda (u), \mu (u))\in
 \mathbb
C\times \mathcal K$, between $\mathcal U^{\mathcal E}$ and $\mathbb C\times \mathcal K ,$ satisfying $$\gamma (u, u')= \overline {\lambda
(u)}+\lambda (u')+\langle \mu (u), \mu (u')\rangle .$$
In the
following, for simplicity of notation, though we have different
product systems,  we will be using same $\lambda $ and $\mu  $ for
the corresponding bijections. This shouldn't cause any confusion. We need couple of lemmas before we state our main theorem. We omit the proof of the first Lemma.

\begin{lem}\label{rank}
Let $\gamma $ be the covariance kernel on the set of all units in a product system $\mathcal E.$ Suppose there is a function $a:\mathcal U^{\mathcal E}\rightarrow \mathbb C$ such that the function $L:\mathcal U^{\mathcal E}\times \mathcal U^{\mathcal E} \rightarrow \mathbb C$ defined by $L(x,y)=\gamma(x,y)-\overline{a(x)}-a(y),$ $x,y \in \mathcal U^{\mathcal E}$ is positive definite, then  $$\sup \{ \mbox{~rank} [ L(x_i,x_j)]_{n\times n}: x_i, x_j \in \mathcal U^{\mathcal E}, n\geq 1 \}= \mbox{index} (\mathcal E)$$
\end{lem}

Recall that units of a spatial product system generate a type I product subsystem. The index of the product system is same as the index of this subsystem. Further, any type I product system is isomorphic to the exponential product system (\cite{Ar01}) or the product system consisting of symmetric Fock spaces $\{ \Gamma (L^2[0,t], K) \} _{t>0}$ where $K$ is a Hilbert space with dim $K$ equal to the index of the product system. In this picture of type I product system, units are parametrized by exponential vectors:
$$\{ e^{qt}e(x\chi|_{t]})\}_{t>0} ~~~(q,x)\in \mathbb C\times K.$$
The automorphisms of this product system is parametrized  by  triples
$\phi :=[q,z,U],$ where $q\in \mathbb R,$ $z\in K$,  $U$ is a unitary in $B(K),$  and $\phi$ acts on the exponential
vectors by
$$\phi e(x\chi|_{t]})=e^{-iqt-\frac{\| z\| ^2t}{2}-\langle z, Ux\rangle t}e((z+Ux)\chi|_{t]}).$$
Then the adjoint of $\phi,$ $\phi^*$ is parameterized by the tuple
$[-q, -U^*z, U^*].$

\begin{lem}\label{indexunit}
Let $A$ be a non-empty subset of a separable Hilbert space $K.$ Then the set of all units in the product subsystem $\mathcal E$ of $\Gamma(L^2[0,t],K),$ generated by units $\{e(x\chi|_{t]}):x\in A\}_{t>0}$ is given by $$\{e^{ct}e(y+x_0)\chi|_{t]}):y\in \overline{\mbox{span}}(A-x_0), c\in \mathbb C\}_{t>0}$$ where $x_0$ is any fixed vector in $A.$ In particular $$\mbox{index}~~(\mathcal E )= \mbox{dim}{~}\overline{\mbox{span}}(A-x_0) $$

\end{lem}
Proof: Fix $x_0\in A.$ As the type I product system is transitive, we
can get an automorphism  $\phi=[q,U,x_0]$ which sends vacuum unit to
$e^{-iqt-\frac{\|x_0\| ^2t}{2}}e(x_0\chi_{t]}).$ Then we claim that it is enough to show the
following assertion : For a subset $B\subseteq K,$ The set
$\{e(y\chi|_{t]}):y \in B\cup {0}\}$  generates units of the form
$$\{e^{\alpha t}e(y\chi|_{t]}):y\in \overline{\mbox{span}}~B,\alpha\in \mathbb C\}.$$
 First assume that the assertion is true. For $x\in A,$
$\phi^*e(x\chi|_{t]})=e^{iqt-\frac{\| x_0 \|^2t}{2}+\langle x_0, x\rangle t}e((U^*(x-x_0))\chi|_{t]}).$ Set
$C=\{U^*(x-x_0):x\in A\}.$ Under $\phi^*,$ The set $\{e(x\chi|_{t])}:
x\in A\}$ maps to $\{e^{iqt-\frac{\| x_0 \|^2t}{2}+\langle x_0, x\rangle t}e(x\chi|_{t])}: x\in C\}.$ As $C=C\cup \{0\},$ by
the assertion we get that the set $\{e(x\chi|_{t]}):x \in C\}$
generates units of the form
$$\{e^{\alpha t}e(x\chi|_{t]}): x \in
\overline{\mbox{span}}~ C,\alpha \in \mathbb C \}.$$ Now under the
image of $\phi,$ we get that the set $\{e(x\chi|_{t]}):x \in A\}$
generates units of the form $$\{e^{\alpha t}e(x\chi|_{t]}): x \in
\overline{\mbox{span}}~ (UC)+x_0,\alpha \in \mathbb C \}.$$ Now ~
$\{x:x\in \overline{\mbox{span}}(UC)+x_0\}=\{y+x_0: y\in \overline{\mbox{span}}(A-x_0)\}.$   Hence the
lemma is proved. Now we will prove the assertion. Let
$K_1=\overline{\mbox{span}}B.$ Now in the subsystem
$\Gamma_{sym}(L^2([0,t],K_1)),$ consider the set
$\{e(z\chi|_{t]}):z\in B\cup \{0 \}\}.$ They generate exponentials of all  step functions taking values in $B\cup \{0 \}.$ As $B$ is total in
$K_1.$ So by a result of Skeide \cite{Sk01}, it is all of
$\Gamma_{sym}(L^2([0,t],K_1)).$ Hence its units are parameterized by
$\{e^{\alpha t}e(y\chi|_{t]}):y\in K_1 , \alpha \in \mathbb C \}.$

Let us denote the product system generated by the exponential vectors $\{e(x\chi|_{t]}): x\in A\}$ by  $\mathcal F$.
We wish to calculate its index. Covariance function is the restriction of the covariance function of $\Gamma_{sym}(L^2([0,t],K_1)),$ which
is $\gamma((\alpha,x),(\beta,y))=\bar{\alpha}+\beta+\langle x,y\rangle.$ Units of $\mathcal F$ are parameterized
by $\mathbb C\times [x_0+\overline{\mbox{span}}(A-x_0)].$ Under the automorphism map $\phi^*,$ $[x_0+\overline{\mbox{span}}(A-x_0)]$ maps to the subspace $U^*(\overline{\mbox{span}}(A-x_0)).$ So we get
$\mbox{index}{~} (\mathcal F) =\mbox{dim}{~}U^*(\overline{\mbox{span}}(A-x_0))=\mbox{dim}{~}(\overline{\mbox{span}}(A-x_0)).$\qed

\begin{lem}\label{ind}
Let $\mathcal E$ be a spatial product system and $Z\subset\mathcal
U^{\mathcal E}$ be a subset of the set of all units in $\mathcal E.$
Let $H_Z$ be the Arveson Hilbert space associated to $Z.$ Then
$\mbox{dim}H_Z=\mbox{ind} ~\mathcal E$ if and only if
$\overline{\mbox{span}}\{u^1_{t_1}\otimes u^2_{t_2}\otimes \cdots  \otimes u^k_{t_k}: 1\leq
i\leq k,u^i\in Z, \sum t_j=t,  k\geq 1\}={\mathcal E}_t$ for all $t>0.$
\end{lem}
Proof: With out loss of generality, we may assume that the given product system
is of type I. We then identify the product system with symmetric Fock space product system.
Then $Z$ can be identified with a subset of $\mathbb C\times
H_{{\mathcal U}^\mathcal E}.$  Take  $A=\{x:(\alpha,x)\in
Z ~~\mbox{for some}~~\alpha \in {\mathbb C} \}$ from the construction of Arveson Hilbert space $H_Z,$ it follows easily that
$H_Z=\overline{\mbox{span}} (A-x_0),$  where $x_0$ is any fixed vector in $A.$ Now the result follows
from Lemma \ref{indexunit}. \qed

\begin{thm} Suppose $\phi =\{ \phi _t: t\geq 0\} $ and $\psi = \{ \psi _t:
t\geq 0\} $ are two $E_0$ semigroups on $\mathcal {B} (\mathcal {H})$ and $\mathcal {B}(\mathcal {K})$
respectively and $U=\{ U_t: t\geq 0\}$ and $V=\{ V_t: t\geq 0\}$ are two strongly continuous semigroups of
contractions which intertwine $\phi _t$ and $\psi _t$ respectively. Consider the CP semigroup $\tau _t$ on
$\mathcal {B}( \mathcal {H}\oplus \mathcal {K})$ defined by
 $\tau _t\left(
 \begin{array}{cc}
 X & Y \\
 Z & W \\
 \end{array}
 \right)=\left(
 \begin{array}{cc}
 \phi _t(X) & U_tYV^*_t \\
 V_tZU^*_t & \psi _t(W) \\
 \end{array}
 \right).$
Let $(\mathcal E, B)$, $(\mathcal F, C)$  and $(\mathcal G,W)$ be the Arveson's product systems associated to
$\phi,$ $\psi$ and $\tau.$
 Then the following holds:
\begin{enumerate}
\item  there exist two units $(u^0)_{t>0}$ and $(v^0)_{t>0}$ of $(\mathcal E,B)$ and $(\mathcal F,C)$ respectively such that  $D:=\{ D_t=|u^0_t\rangle\langle v^0_t|:t>0\}$ from $\mathcal F$ to
$\mathcal E$ such that $\mathcal E\otimes_D\mathcal F=\mathcal G.$

\item The type I part of the amalgamated product is the amalgamated product of
type I parts of $\mathcal E$ and $\mathcal F$: $(\mathcal
{E}^I\otimes _D\mathcal {F}^I)=(\mathcal E\otimes _D\mathcal F)^I;$
\item $$\mbox{index}~ (\mathcal E\otimes _D\mathcal F)= \left\{ \begin{array}{ll}
\mbox{index}~{(\mathcal E)}+ \mbox {index} ~(\mathcal F) & ~~\mbox {if}~~\|u^0_t\|=\|v^0_t\|=1 \forall t>0; \\
\mbox {index}~(\mathcal E)+\mbox{index} ~(\mathcal F)+1
&~~\mbox{otherwise.}
\end{array}\right.
$$
\end{enumerate}
\end{thm}

Proof:  Strong continuity properties of $U$ and $V$ imply that the
product system associated to the $E_0$ dilation of $\tau_t$ is an
Arveson's product system. Let $(G,\delta)$ be the inclusion system of
$\tau.$ Now from  Theorem \ref{amalgamation}, $G_t=\mathcal
E_t\oplus_{D_t} \mathcal F_t.$ We conclude that $\mathcal G=\mathcal
E\otimes_D\mathcal F,$ where
$$D_t=P_{\mathcal E_t}\pi _t(|a\rangle\langle b|)|_{\mathcal F_t}=|U_ta\rangle\langle V_tb|,$$ (Here $\pi _t$ denotes the minimal dilation of $\tau _t.$) Take
$u^0_t=U_ta$ and $v^0_t=V_tb.$ Then clearly $(u^0_t)_{t>0}$ and
$(v^0_t)_{t>0}$ are units of $(\mathcal E,B)$ and $(\mathcal F,C)$
respectively.


 Consider $M_t =\overline{\mbox{span}}\{\left[
                          \begin{array}{c}
                            u_t \\
                            0 \\
                          \end{array}
                        \right],\left[
                                  \begin{array}{c}
                                    0 \\
                                    v_t \\
                                  \end{array}
                                \right], u \in \mathcal U^\mathcal E , v \in \mathcal U^\mathcal F\}
.$ As $u \in \mathcal U^\mathcal E$ and $v \in \mathcal U^\mathcal
F$, $\left[ \begin{array}{c}
                                                                          u \\
                                                                          0 \\
                                                                        \end{array}
                                                                      \right]$ and $\left[
                                                                                    \begin{array}{c}
                                                                                      0 \\
                                                                                      v \\
                                                                                    \end{array}
                                                                                  \right]$
 are strong units in $ (G, \delta ),$ and  $(M, \delta|_M)$ is an inclusion subsystem of $(G,\delta)$. Let  $(\mathcal
M,W|_\mathcal M)$ be its generated product system, we get that
$(\mathcal M,W|_\mathcal M)\supset\mathcal E^I\otimes_D\mathcal
F^I$. As $(\mathcal M ,W|_\mathcal M)$ is a product subsystem of
$(\mathcal G^I,W),$ we have $\mathcal E^I\otimes_D\mathcal F^I$ as a
product subsystem of $(\mathcal G^I,W)$. In particular amalgamation
of two type I system is again type I. Now let $\left[
                                                                                                \begin{array}{c}
                                                                                                  u \\
                                                                                                  v \\
                                                                                                \end{array}
                                                                                              \right]
$ be a unit of $(G,\delta)$. $D_s=|u^0_s\rangle\langle v^0_s|,$
implies $\tilde{D_t}$ is closed. So invoking the Lemma \ref{formula}
we can conclude that $\left[
                                                                                                \begin{array}{c}
                                                                                                  u_t \\
                                                                                                  v_t\
                                                                                                \end{array}
                                                                                              \right]
\in M_t$ for every $t>0.$ In particular it implies that $(G^I,W)$ is
a product subsystem of $(\mathcal M, W|_\mathcal M).$ So $(\mathcal
E^I\otimes_{D}\mathcal F^I)\simeq\mathcal G^I$, proving (2).


So from the  argument above, it follows that the set
$$Z=\{ \widehat{\left[
\begin{array}{c}u\\ 0\end{array}\right]} : u\in \mathcal U^{\mathcal
E }\}\bigcup \{\widehat{ \left[
\begin{array}{c}0\\ v\end{array}\right]} : v\in \mathcal U^{\mathcal
F}\}\subset \mathcal U^{\mathcal G}$$ generates $\mathcal G^I$ (Here  the lift to the generated product system of a unit $x$ of an inclusion system is denoted by $\hat {x}.$) So
from \ref{ind}, $\mbox{ind}(\mathcal G)=\mbox{dim}H_Z.$ So it is
enough to calculate the rank of the covariance kernel on $Z.$ The covariance function can be computed as follows. For
arbitrary units $u\in \mathcal U^{\mathcal E}, v\in \mathcal
U^{\mathcal F},$
$$\langle \widehat{\left[
\begin{array}{c}u_t\\ 0\end{array}\right]}, \widehat{\left[
\begin{array}{c}0\\ v_t\end{array}\right]}
\rangle = \langle u_t, |u^0_t\rangle \langle v^0_t|v_t\rangle =
e^{t(\gamma (u, u^0)+ \gamma (v^0, v))}.$$ Therefore,
$$\gamma (\widehat{\left[
\begin{array}{c}u\\ 0\end{array}\right]}, \widehat{\left[
\begin{array}{c}0\\ v\end{array}\right]}\rangle =
\overline {\lambda (u)}+\lambda (u^0)+\langle \mu  (u), \mu
(u^0)\rangle+\overline {\lambda (v)}+\lambda (v)+\langle \mu (v^0),
\mu  (v)\rangle.$$ Similarly, for $u, u'\in \mathcal U^{\mathcal E}$
and $v, v'\in \mathcal U^{\mathcal F},$
 \Bea \gamma (\widehat{\left[
\begin{array}{c}u\\ 0\end{array}\right]}, \widehat{\left[
\begin{array}{c}u'\\ 0\end{array}\right]}\rangle &=&
\overline {\lambda (u)}+\lambda (u')+\langle \mu  (u), \mu
(u')\rangle , \\
\gamma (\widehat{\left[
\begin{array}{c}u\\ 0\end{array}\right]}, \widehat{\left[
\begin{array}{c}0\\ v\end{array}\right]}\rangle &=&
\overline {\lambda (v)}+\lambda (v')+\langle \mu (v), \mu
(v')\rangle.\Eea Now take $Y=\mathcal U^{\mathcal E}\bigcup \mathcal
U^{\mathcal F}$ and define $a:Y\to \mathbb C$ by $$ a(u)= \lambda
(u)+\langle \mu  (u^0), \mu  (u)-\frac{1}{2}\mu  (u^0)\rangle
~~\mbox {if} ~~u\in \mathcal U^{\mathcal E}$$ and
$$a(v)=  \lambda (u^0)+\frac{1}{2}\langle \mu  (u^0), \mu  (u^0)\rangle + \overline{\lambda
(v^0)}+\lambda (v)+\langle \mu  (v^0), \mu  (v)\rangle ~~\mbox {if}
~~v\in \mathcal U^{\mathcal F}.$$ Similarly define $L:Y\times Y\to
\mathbb C$ by $$L(x, y)= \gamma (\widehat{[x]}, \widehat{[y]})-\overline
{a(x)}-a(y).$$ (Here $[x]$ denotes $\left[
\begin{array}{c}x\\ 0\end{array}\right]$ or $\left[
\begin{array}{c}0\\ x\end{array}\right]$ according as $x$ is in
$\mathcal U^\mathcal E$ or $\mathcal U^\mathcal F.$)  Then by direct
computation: For $u, u'\in \mathcal U^{\mathcal E}$ and $v, v'\in
\mathcal U^{\mathcal F},$ \Bea L(u,u')&=& \langle \mu  (u)-\mu
(u^0), \mu  (u')-\mu
(u^0)\rangle \\
L(u, v)&=& 0\\
L(v, v') &=&  \langle \mu  (v)-\mu  (v^0), \mu  (v')-\mu
(v^0)\rangle +p \Eea where $p:=-[\overline {\lambda
(u^0)}+\lambda(u^0)+\overline {\lambda(v^0)}+\lambda(v^0)+ \langle
\mu  (u^0), \mu  (u^0)\rangle + \langle \mu  (v^0), \mu
(v^0)\rangle]= -[\gamma (u^0, u^0)+\gamma (v^0, v^0)]$. It is to be
noted that  as $e^{-tp}= \| u^0_t\|.\|v^0_t\| \leq 1 $ for all $t$,
$p$ is non-negative and $p=0$ iff $\|u^0_t\|=\|v^0_t\|=1$ for all
$t.$ So taking direct sum of the range space of $\mu $ with $\mathbb C ,$ we get
$$L(v, v')= \langle (\mu  (v)-\mu  (v^0))\oplus \sqrt {p}, (\mu  (v')-\mu
(v^0))\oplus \sqrt (p)\rangle .$$

For any unit $u\in \mathcal U^{\mathcal E}$, we can find another
unit $\tilde u \in  \mathcal U^{\mathcal E}$ such that $\mu (u)=\mu
(\tilde u)-\mu  (u^0).$ Then it is clear that maximal rank of
$[L(x_i, x_j)]$ with $x_1, \ldots , x_n \in \mathcal U^{\mathcal E}$
is equal to index $(\mathcal E)$. Similarly, maximal rank of
$[L(y_i, y_j)]$ with $y_1, \ldots , y_n \in \mathcal U^{\mathcal F}$
is equal to index $(\mathcal F)+1$, if $p>0$ and is equal to index
$(\mathcal F)$ if $p=0.$ The theorem follows from the Lemma
\ref{rank}.\qed

{\bf Acknowledgements:} B. V. Rajarama Bhat thanks UKIERI for
financial support.


\bibliographystyle{amsplain}

\end{document}